\documentclass[12pt]{article}
\usepackage{amsmath,amsthm,amsfonts,graphicx,amssymb}
\usepackage{url}
\newtheorem{Theorem}{Theorem}[section]
\newtheorem{lemma}[Theorem]{Lemma}
\newtheorem{corollary}[Theorem]{Corollary}

\theoremstyle{definition}
\newtheorem{definition}[Theorem]{Definition}

\numberwithin{equation}{section}

\def\bp{\noindent{\bf Proof.}  \ }
\newcommand{\ep}{\hfill $\square$}

\begin{document}

\title{Laplacian Spectra of  Regular Graph Transformations
\thanks{Aiping Deng and   Juan Meng are supported in part by the Fundamental Research Funds for the Central Universities of China 11D10902 and
11D10913. }}

\author{ Aiping Deng$^{\rm a}$\footnote{Corresponding author.
Email: apdeng@dhu.edu.cn. Tel: 86-21-67792089-568. Fax: 86-21-67792311 }\ ,
 Alexander Kelmans$^{\rm b,c}$,  Juan Meng$^{\rm a}$ \\
$^{\rm a}${\small\em{Department of Applied Mathematics, Donghua University, 201620 Shanghai, China} } \\
$^{\rm b}${\small\em{Department of Mathematics, University of Puerto Rico, San Juan, PR, United States}}\\
$^{\rm c}${\small\em{Department of Mathematics, Rutgers University, New Brunswick, NJ, United States}}
  }

\date{}
\maketitle

\begin{abstract}
Given a graph $G$ with vertex set $V(G)= V$ and edge set $E(G) = E$, let $G^l$ be the line graph and $G^c$ the complement of $G$.
Let $G^0$ be the graph with $V(G^0) = V$ and with no edges, $G^1$ the complete graph with the vertex set $V$, $G^+ = G$ and
$G^- = G^c$.
Let $B(G)$ ($B^c(G)$) be the graph with the vertex set $V\cup E$ and
such that $(v,e)$ is an edge in
$B(G)$ (resp., in $B^c(G)$)
if and only if $v \in V$, $e \in E$ and  vertex $v$ is incident (resp., not incident) to edge $e$ in $G$.  Given $x, y, z \in \{0,1, +, -\}$,
the {\em ${xyz}$-transformation $G^{xyz}$  of}
$G$
is the graph with the vertex set $V(G^{xyz}) = V \cup E$ and the edge set $E(G^{xyz}) = E(G^x) \cup E((G^l)^y) \cup E(W)$, where
 $W = B(G)$  if $z = +$,
 $W = B^c(G)$
 if $z = -$, $W$ is the  graph with
 $V(W) = V \cup E$ and with no edges  if $z=0$, and
$W$ is the complete bipartite graph with parts $V$ and $E$ if $z = 1$.
In this paper we obtain the Laplacian characteristic polynomials and some other Laplacian parameters of
every $xyz$-transformation of an $r$-regular graph $G$ in terms of $|V|$,  $r$, and the Laplacian spectrum
 of $G$.
\\[1.5ex]
\\[0.7ex]

\end{abstract}

\section{Introduction}

\indent

The graphs in this paper are simple and undirected.
All notions on graphs and matrices that are used but not defined here can be found in
\cite{B&M07,D,G,horn,west}.
\\[1ex]
\indent
Let ${\cal G}$ denote the set of simple undirected graphs.
Various important results in graph theory have been obtained by considering some functions
$F: {\cal G} \to  {\cal G}$
or  $F_s: {\cal G}_1\times  \ldots \times {\cal G}_s \to  {\cal G}$ called
{\em operations} or {\em transformations}
(here  each ${\cal G}_i = {\cal G}$) and
by establishing how these operations affect certain properties or parameters of graphs.
The complement,  the $k$-th power of a graph, and  the line graph are well known  examples of such operations.
The  Bondy-Chv\'atal and Ryz\'{a}\u{c}ek closers of graphs are very useful operations in graph Hamiltonicity theory
\cite{B&M07}.
(Strengthenings and extensions of the Ryz\'{a}\u{c}ek result are given in \cite{Kclosure}).
Some graph operations introduced by A. Kelmans
(see, in particular,  \cite{Ktrlp, KoperR}) turn out to be
monotone with respect to various partial order relations on the set of graphs.
For that reason these operations
turned out to be very  useful in obtaining
non-trivial results on graphs  of given size with various  extreme properties (with the  maximum number of spanning trees and some other Laplacian parameters of graphs, with the maximum reliability of graphs having randomly deleted edges, etc.),
see, for example,   \cite{Kmxtr, Kprobcmpr}.
The operation of voltage lifting on a base graph introduced by Gross and Tucker can be generalized to digraphs
\cite{D&W05, G&T87}.
Using this operation one can obtain the derived covering (di)graph and deduce the relationship between the adjacency characteristic polynomials of the base (di)graph and its derived covering (di)graph \cite{D&W05, DS&W07, M&S95, wuyk}.
\\[1ex]
\indent
In this paper we consider  certain graph operations
depending on parameters  $x, y, z \in \{0,1,+,-\}$.
These operations  induce  functions $T^{xyz}: {\cal G} \to {\cal G}$. We put  $T^{xyz}(G) = G^{xyz}$ and call
$G^{xyz}$ the {\em $xyz$-transformation of $G$}.
We describe for all  $x, y, z \in \{0,1,+,-\}$
the  Laplacian characteristic polynomials and some other Laplacian parameters of
$xyz$-transformations of an $r$-regular graph $G$.
This descriptions revealed the following fact interesting in itself: if $G$ is $r$-regular, then
the Laplacian spectrum of  $G^{xyz}$ is uniquely defined by $|V(G)|$, $r$,
and the Laplacian spectrum of $G$; moreover, the Laplacian eigenvalues are the roots of a quadratic polynomial with the coefficients depending on $|V(G)|$, $r$, and the Laplacian spectrum of $G$.
Furthermore, for $(xyz) \in \{(00+),(0++),(+0+)\}$
the number of spanning trees of $G^{xyz}$ are uniquely defined by  $|V(G)|$, $r$, and the number of spanning trees of $G$ (see Theorem \ref {trees} and Corollaries \ref{t(G00+)}, \ref{t(G+0+)}, and \ref{t(G0++)} below).
The approach we have used to obtain all these formulas  may also be useful in further research along this line.  The results of this paper may be considered as a natural and useful extension of Section 2  ``Operations on Graphs and the Resulting  Spectra'' in book \cite{CDS}.

The Reciprocity Theorem \cite{Ktree1} (see also Theorem \ref{reciprocity} below) provides for every graph $G$ the relation between the Laplacian characteristic polynomial of $G$ and its complement $G^c$. For that reason it is sufficient to describe the  Laplacian characteristic polynomials of graph $xyz$-transformations up to
the graph operation of taking the complement.
\\[1ex]
\indent
In Section \ref{notions} we introduce main
 notions, notation, simple observations, and  some preliminaries.
In Section \ref{0-in-xyz} we describe (up to the complementarity) the Laplacian characteristic polynomials  of the transformations $G^{xyz}$ with $\{x,y,z\}\cap \{0, 1\} \ne \emptyset $.
In Section \ref{x,y,z-in-(+,-)} we  describe
the Laplacian characteristic polynomials  of transformations $G^{xyz}$ with $\{x,y,z\}\cap \{0, 1\} = \emptyset $, i.e. with  $x, y, z \in \{+, -\}$.
In Sections \ref{examples} we consider some  transformations of cycles and show that different transformations of the same graph may be isomorphic. Section \ref{remarks} contains some additional remarks.
In the Appendix we provide for all $x, y, z \in \{0,1,+, -\}$ the list of formulas
for the Laplacian characteristic polynomials and the number of spanning trees of the $(xyz)$-transformations of an $r$-regular graph $G$ in terms of $r$,
the number of vertices, the number of edges, and the Laplacian spectrum of $G$. This catalog may be pretty useful for further research in the spectrum graph theory.

\section{Some notions, notation, and  preliminaries}
\label{notions}

\indent

Let $G =(V, E)$ be a graph with vertex set
$V = V(G)$ and edge set
$E = E(G)$. Let $v(G) = |V(G)|$ and $e(G) = |E(G)|$.
 The {\em degree} $d(v, G)$ {\em of vertex} $v \in V$  is the number of vertices in $G$ adjacent to $v$.
 Let $t(G)$ denote the number of spanning trees of $G$.

Given two graphs $G$ and $H$, an  {\em isomorphism from $G$ to $H$} is a bijection $\alpha $ from $V(G)$ to $V(H)$ such that $(u, v) \in E(G) \Leftrightarrow (\alpha (u), \alpha (v)) \in E(H)$.

 The {\em complement} $G^c$ of a graph $G$ is the graph with vertex set $V(G^c) = V(G)$ and $(u, v) \in E(G^c) \Leftrightarrow (u, v) \not \in E(G)$ for any
 $u, v \in V(G)$ and $u \ne v$.

\begin{figure}[ht]
\begin{center}
\scalebox{0.7}[.7]{\includegraphics{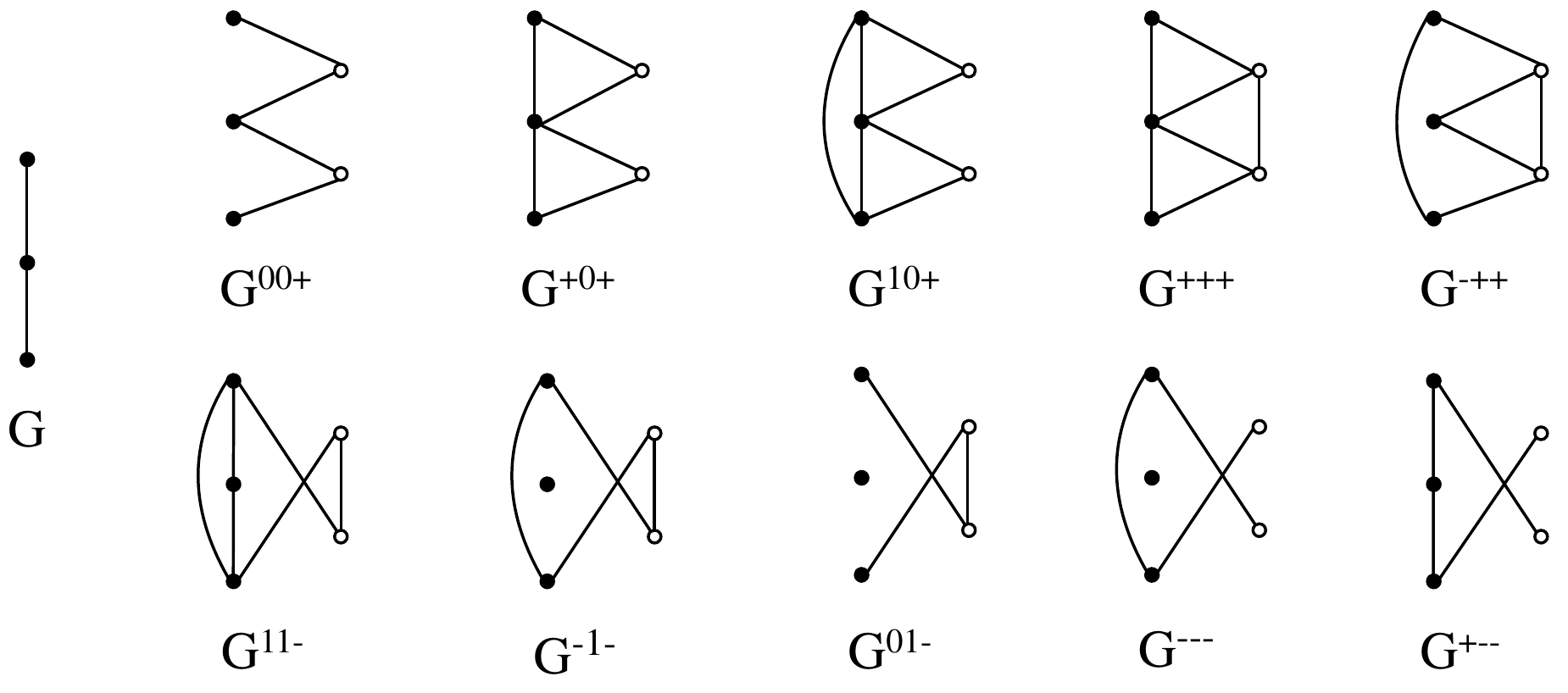}}
\end{center}
\caption{Some $xyz$-transformations of a 3-vertex path.}\label{tran}
\end{figure}

 The {\em line graph} $G^l$ of a graph $G$ is the graph
with vertex set $E(G)$ and two vertices are adjacent in $G^l$ if and
only if the corresponding edges in $G$ are adjacent.
%
%
Graphs $G$ and $H$  are called {\em isomorphic} if there exists an isomorphism from $G$ to $H$.
\\[1ex]
\indent
For a graph
$G = (V, E)$, let
$G^0$ be the graph with $V(G^0) = V$ and  with no edges, $G^1$ the complete graph with $V(G^1) = V$,
$G^+ = G$, and $G^- = G^c$.
Let $B(G)$ ($B^c(G)$) be the graph with the vertex set $V\cup E$ and
such that $(v,e)$ is an edge in
$B(G)$ (resp., in $B^c(G)$) if and only if $v \in V$, $e \in E$, and  vertex $v$ is incident (resp., not incident) to edge $e$ in $G$.
 For example, in Figure \ref{tran},
$G^{00+} = B(G)$ and $B^c(G)$ is obtained from
$G^{01-}$ by deleting  the edge connecting two white vertices.
\\[1ex]
\indent
The graph transformations we are going to discuss are defined as follows.
\begin{definition}
Given a graph $G = (V, E)$ and three variables
$x, y, z \in \{0,1, +, -\}$, the {\em ${xyz}$-transformation $G^{xyz}$  of}
$G$ is the graph with the vertex set
$V(G^{xyz}) = V \cup E$
and the edge set $E(G^{xyz}) = E(G^x) \cup E((G^l)^y) \cup E(W)$, where
 $W = B(G)$  if $z = +$, $W = B^c(G)$ if $z = -$, $W$ is the  graph with $V(W) = V \cup E$ and with no edges if $z=0$, and
$W$ is the complete bipartite graph with parts $V$ and $E$ if $z = 1$.
\end{definition}

Examples of ${xyz}$-transformations of a 3-vertex path  are given in Figure \ref{tran}.
Graphs $G^{+++}$ and $G^{00+}$ are  called in \cite{CDS} the {\em total graph} and the {\em subdivision graph} of $G$, respectively.
\\[1ex]
\indent
Let $G$ be a graph with vertex set
$V  = \{v_1, \ldots , v_n\}$ and edge set
$E   = \{e_1, \ldots ,e_m\}$.
  The {\em incidence matrix} $Q(G)$ of $G$ is
 the $(V \times E)$-matrix $\{q_{ij}\}$, where
$q_{ij} = 1$ if vertex $v_i$ is incident to edge $e_j$ and
$q_{ij} = 0$, otherwise.
Let $A(G)$ be the $(V\times V)$-matrix $\{a_{ij}\}$, where
$a_{ij} = 1$ if $(v_i, v_j) \in E$ and $a_{ij} = 0$, otherwise.
Let $D(G)$ be the (diagonal)  $(V\times V)$-matrix
$\{d_{ij}\}$, where $d_{ii} = d(v_i,G)$ and $d_{ij} = 0$ for $i \ne j$.
The matrices $A(G) $, $D(G)$ and $L(G) = D(G) - A(G)$ are called
the {\em adjacency matrix}, the {\em degree matrix}, and the {\em Laplacian matrix} of $G$, respectively.
The {\em adjacency polynomial},
the {\em adjacency spectrum} and the
{\em adjacency eigenvalues} of $G$ are
the characteristic polynomial
$A(\alpha, G) = \det (\alpha I - A(G))$,
the spectrum, and the eigenvalues
of $A(G)$, respectively.
Similarly, the
{\em Laplacian polynomial}, the
{\em Laplacian spectrum} and the
{\em Laplacian eigenvalues} of $G$ are
the characteristic polynomial
$L(\lambda , G) = \det (\lambda I - L(G))$,
the spectrum, and
 the eigenvalues of $L(G)$, respectively.
Let $I_n$ be the identity $(n\times n)$-matrix and
$J_{m n}$ the all-ones $m\times n$-matrix.

Since $A(G)$ and $L(G)$ are symmetric matrices,
their eigenvalues are real numbers.
It is easy to see that each row sum of $L(G)$  is equal to zero and that
$L(G)$ is a symmetric positive semi-definite matrix
\cite{CDS,kelmans}. Therefore all eigenvalues $\lambda _i(G)$  of $L(G)$
are  real non-negative numbers and one of them is equal to zero; we order them in the descendant order:
\[
\lambda _1(G) \geq \lambda_2(G) \geq \cdots \geq \lambda_n(G) = 0.\]
The set  $Sp(G) = \{\lambda _1(G), \ldots ,  \lambda_n(G) \}$ is called the {\em Laplacian spectrum of $G$}.
%
\\[1ex]
\indent
Some graph properties
of the transformations $G^{xyz}$ with $x,y,z\in \{+,-\}$ have been discussed and obtained in \cite{linshu, wubaoy2, wumeng1}.
\\[0.7ex]
\indent
For a regular graph $G$, the adjacency characteristic polynomials  and the adjacency spectrum of $G^{00+}$,
$G^{+0+}$, $G^{0++}$ and the total graph $G^{+++}$ are given in \cite{CDS} (pages 63 and 64).
The adjacency characteristic polynomials  and the adjacency spectrum  of the other seven $G^{xyz}$ with $x, y, z \in \{ +, -\}$ are obtained in
\cite{yanxu}.
The definition of $xyz$-transformation can be easily extended to digraphs, which is also a generalization of the digraph transformations defined by Liu and Meng \cite{liumeng}. Zhang, Lin, and Meng have described
the adjacency characteristic polynomials of $D^{00+}, D^{+0+}, D^{0++}$ and the total digraph $D^{+++}$  for any digraph $D$ \cite{zhang}. The adjacency characteristic polynomials of other $D^{xyz}$ of a regular digraph $D$ with $x, y, z \in \{+,-\}$ are obtained in \cite{liumeng}.
\\[0.7ex]
\indent
Very few results are known for the Laplacian spectra of transformations.
In 1967 A. Kelmans published the following results
on the Laplacian polynomial of $G^{0++}$, $G^{0+0}$,
$G^{00+}$, and $G^l$ for  a regular graph $G$.
These results are included in the survey papers
 \cite{mohar} (Theorem 3.8) and \cite{newman} (Theorem 1.4.2) with an error, namely,
graph $G^{0++}$ is mistakenly called the total graph of
$G$.
\begin{Theorem}
\label{G0++}
{\em \cite{kelmans}}
Let $G$ be an $r$-regular graph with $n$ vertices and $m$ edges. Then
\[L(\lambda , G^{0++}) =
(\lambda -  r - 1)^n(\lambda - 2r - 2)^{m - n}
L(\frac{\lambda ^2 -  (r + 2)\lambda}{\lambda- r - 1}, G).\]
\end{Theorem}
\begin{Theorem}
\label{Ln(G),G0+0}
{\em \cite{kelmans}}
Let $G$ be an $r$-regular graph with $n$ vertices and $m$
edges. Then
\[L(\lambda , G^l) =
(\lambda - 2r)^{m - n} L(\lambda , G)~~and~~
L({\lambda , G^{0+0}}) =  \lambda ^n L(\lambda , G^l).\]
\end{Theorem}
\begin{Theorem}
\label{G00+}
{\em \cite{kelmans}}
Let $G$ be an $r$-regular graph with $n$ vertices and $m$ edges. Then
\[L(\lambda , G^{00+}) = (-1)^n(\lambda - 2)^{m - n}
L(\lambda (r + 2 - \lambda), G).\]
\end{Theorem}

Since $nt(G) = (-1)^{n-1}\lambda ^{-1}L(\lambda , G)|_{\lambda = 0}$, where $n = v(G)$ \cite{Ktree1},
we have  from Theorems \ref{G0++}, \ref{Ln(G),G0+0}, and \ref{G00+}:
\begin{Theorem}
\label{trees}
{\em \cite{kelmans}}
Let $G$ be an $r$-regular graph with $n$ vertices and $m$
edges. Then
\[t (G^{0++}) =
\frac{n}{m+n} 2^{m - n}(r + 1)^{m-1}(r+2)   t(G),\]
\[ t(G^{00+}) = \frac{n}{m+n}2^{m - n} (r+2)  t( G),\]
and
\[ t(G^l) = \frac{n}{m} 2^{m - n} r^ {m - n}  t(G).\]
\end{Theorem}

A set ${\cal A}$ of graphs is {\em closed under complementarity} if for every graph in ${\cal A}$ its complement is also in ${\cal A}$.
Given a graph $G$ and $S \subseteq \{0,1, +, -\}$, let
${\cal F}(G, S)$ denote the set of graphs $G^{xyz}$ such that $x, y, z \in S$.
If ${\cal F}(G, S)$ is closed under complementarity, then in order to find the Laplacian polynomial for the graphs in
${\cal F}(G, S)$, it is sufficient to find the solutions for a ``half'' of graphs in ${\cal F}(G, S)$ and to obtain the solutions for the graphs of the other ``half'' using the  {\em Reciprocity Theorem} \ref{reciprocity} \cite{Ktree1}.
It is easy to see that  ${\cal F}(G, \{+, -\})$ and ${\cal F}(G, \{0, 1, +, -\})$ (and therefore ${\cal F}(G, \{0, 1, +, -\} \setminus \{+, -\})$) are closed under complementarity.
Hence the Reciprocity Theorem below can be used for these
classes of transformations.

\begin{Theorem}
{\em \cite{Ktree1}}
\label{reciprocity}
Let $G$ be a simple graph with $n$ vertices.
Then
\\[1ex]
$(a1)$ $\lambda _i(G) + \lambda _{n - i}(G^c) = n$
for every $i \in \{1, \ldots , n-1\}$ or, equivalently,
\\[1ex]
$(a2)$
$(n - \lambda ) L(\lambda , G^c) =
(-1)^{n-1} \lambda  L(n - \lambda , G)$.

\noindent Moreover, the matrices $L(G)$ and $L(G^c)$ are simultaneously diagonalizable.
\end{Theorem}
Similar results for weighted graphs are obtained in \cite{KeigvMtr}.
\\[1ex]
\indent
We will also use the following two classical and simple facts on matrices.
\begin{lemma}
\label{lemABCD}
{\em \cite{G,horn}}
Let $A$ and $D$ be square matrices. Then
 \[ \left|\begin{array}{cc}
A &\quad B\\ C &\quad D \end{array}\right|\quad
=\left\{\begin{array}{ll}|A|~|D - C A^{-1} B|, &\quad \,\, if\, $A$
\mbox{\, is\, invertible},
\\[1ex]
|D|~|A - B D^{-1}C|, &\quad \,\, if\, $D$
\mbox{\, is\, invertible}.
\end{array}\right.\]
\end{lemma}

The following two useful and very well known lemmas are obvious.
\begin{lemma}
\label{lemQ,Qt}
Given a graph $G$ with $m$ edges, let
$G^l$ be the line graph of $G$ and $Q$ the incidence matrix of $G$. Then
\\[1ex]
$(a1)$
$QQ^{\top} = D(G) + A(G)$ and
\\[1ex]
$(a2)$
$Q^{\top}Q = 2I_m + A(G^l)$.
\end{lemma}

\begin{lemma}
\label{lemJ}
Let $G$ be an $r$-regular  graph with $n$ vertices and $m$ edges and let $A$ and $Q$ be the adjacency and the incidence matrix of $G$, respectively.
Let $k$ be a positive integer. Then
\\[1ex]
$(a1)$
$Q^{\top}J_{n k} = 2J_{m k}$,
\\[1ex]
$(a2)$
$QJ_{m k}=r J_{n k}$,
\\[1ex]
$(a3)$
$J_{k m}Q^{\top}=rJ_{k n}$,
\\[1ex]
$(a4)$
$J_{k n}Q = 2J_{k m}$,
\\[1ex]
$(a5)$
$J_{kn}A = r J_{k n}$, and
\\[1ex]
$(a6)$
$A J_{n k} = r J_{n k}$.
\end{lemma}

\begin{lemma}\label{Corfg,AJ}
Let $G$ be an $r$-regular graph with $n$ vertices, $A(G) = A$, and
\\
$\lambda_1 \geq \lambda_2 \geq \cdots \geq \lambda_n =0$ the list of the Laplacian eigenvalues of $G$.
Let $P(x,y)$  be a polynomial with two variables and real coefficients.
 Then matrix $P(A, J_{n n})$ has the eigenvalues $\sigma _n = P(r, n)$ and $\sigma _i = P(r - \lambda_i, 0)$ for
 $i=1,2,\cdots, n-1.$
\end{lemma}

\bp ~
Since the Laplacian matrix $L = L(G)$ is symmetric and real, there is a list $B = \{X_1, X_2, \cdots, X_n\}$ of mutually orthogonal  eigenvectors of $L$, where $X_i$ corresponds to $\lambda_i$ and $X_n = J_{n 1}$. Since
$A = rI_n - L$, clearly $AX_i= (r - \lambda_i)X_i$ for each $i$. Since $B$ is an orthogonal basis,
$J_{n n}X_i = 0$ for each $i\neq n$.
Clearly, $J_{n n}J_{n 1} = n J_{n1}$ and
$J_{n n}^2 = n J_{n n}$. By Lemma \ref{lemJ},
$A J_{n n} = J_{n n} A = r J_{n n}$.
Therefore,   $P(A,J_{n n})~X_n = P(r, n)~ X_n$ and
$P(A, J_{n n})~X_i = P(r - \lambda _i, 0)~X_i$ for each $i \ne n$.
\ep
\\[1ex]
\indent
The arguments in this proof are similar to those in
\cite{KeigvMtr}.

\section{Laplacian spectra of   $G^{xyz}$ with
$\{x,y,z\}\cap \{0, 1\} \ne \emptyset $}
\label{0-in-xyz}

\indent

Given a graph $G$ with   $n$ vertices and $m$ edges, we always denote by $A$, $D$ and $Q$, the adjacency matrix, the degree matrix and
the incidence matrix of $G$, respectively, and so if
$G$ is an $r$-regular graph, then
$D = rI_n$ and $2m = rn$.
We put $\lambda _i(G) = \lambda _i$, and so
$\lambda_1\ge \lambda_2 \ge \cdots \ge  \lambda_n = 0$ is the list of the Laplacian eigenvalues of $G$.

\subsection{Laplacian spectra of $G^{xyz}$ with $z = 0$}

\indent

We start with the following simple observation.
\begin{Theorem}
\label{xy0general}
Let $G$ be a graph with $n$ vertices and $m$ edges and
let $x,y \in \{0,1,+,-\}$. Then
$L(\lambda ,G^{xy0}) =
L(\lambda ,G^x)L(\lambda ,(G^l)^y)$.
\end{Theorem}

Since $L(\lambda, G^0) = \lambda^n$, $L(\lambda, G^{+}) = L(\lambda, G)$ and $L(\lambda, G^{1})= \lambda (\lambda -n)^{n-1}$,
we can calculate $L(\lambda, G^{xy0})$ for $x, y \in \{0, 1, +, -\}$ from
Theorems \ref{Ln(G),G0+0},
\ref{reciprocity}  and \ref{xy0general}.
\begin{Theorem}
\label{xy0}
Let $G$ be an $r$-regular  graph with $n$ vertices and $m$ edges.
  Then
\\[1.5ex]
$(a1)$
$L(\lambda, G^{x00}) = \lambda^m L(\lambda, G^x)$ and $L(\lambda, G^{x10}) =
\lambda (\lambda - m)^{m-1} L(\lambda, G^x)$
for $x \in \{0, 1, +\}$,
\\[1.5ex]
$(a2)$
$L(\lambda, G^{-00}) = (-1)^n
(\lambda - n)^{-1} \lambda^{m+1} L(n-\lambda, G)$,
\\[1.5ex]
$(a3)$
$L(\lambda, G^{-10}) = (-1)^n \lambda^2 (\lambda -n)^{-1}(\lambda -m)^{m-1} L(n-\lambda, G)$,
\\[1.5ex]
$(a4)$
$L(\lambda, G^{0-0}) = (-1)^n (\lambda - m)^{-1} \lambda^{n+1} (\lambda - m + 2r)^{m-n} L(m-\lambda, G)$,

\hskip 0.3cm $L(\lambda, G^{1-0}) = (-1)^n (\lambda - m)^{-1} \lambda^2 (\lambda -n)^{n-1} (\lambda - m + 2r)^{m-n} L(m-\lambda, G)$,
\\[1.5ex]
$(a5)$
$L(\lambda ,G^{++0}) =
(\lambda - 2r)^{m - n} L(\lambda , G)^2$,

\hskip 0.3cm $L(\lambda ,G^{-+0}) =
(-1)^n (\lambda - 2r)^{m - n}
\lambda (\lambda - n)^{-1}
L(n - \lambda , G) L(\lambda , G) $,
\\[1.5ex]
$(a6)$
$L({\lambda , G^{0+0}}) =  \lambda ^n (\lambda - 2r)^{m - n} L(\lambda , G)$,

\hskip 0.3cm $L(\lambda ,G^{1+0}) =
\lambda (\lambda - n)^{n-1}(\lambda - 2r)^{m - n} L(\lambda , G)$,
\\[1.5ex]
$(a7)$
$L(\lambda ,G^{+-0}) =  (-1)^n
\lambda (\lambda - m)^{-1}
(\lambda - m + 2r)^{m - n} L(m - \lambda , G)
L(\lambda , G)$ and

\hskip 0.3cm $L(\lambda ,G^{--0}) =
\lambda ^2( \lambda - m)^{-1}(\lambda - n)^{-1}
(\lambda - m +  2r)^{m - n}
L(m - \lambda , G) L(n - \lambda , G)$.
\end{Theorem}

\subsection{Laplacian spectra of $G^{xyz}$ for $z = +$ and $x, y \in \{0, 1\}$}

\begin{Theorem}
\label{00+}
Let $G$ be an $r$-regular graph with $n$ vertices and $m$ edges. Then
\[
L (\lambda ,G^{00+}) =
\lambda (\lambda -r - 2)(\lambda - 2)^{m-n}
\prod_{i=1}^{n-1}\{\lambda ^2-\lambda (r+2)+\lambda _i\}
\]
or equivalently,
\[L(\lambda , G^{00+}) = (-1)^n(\lambda - 2)^{m - n}
L(\lambda (r + 2 - \lambda)), G).\]
\end{Theorem}
\bp
It is easy to see that
\[ A(G^{00+})= \left( \begin{array}{cc} 0 &\quad
Q
\\[1ex] Q^{\top}  &\quad 0
\end{array} \right) \hbox{\ and\ } D(G^{00+})= \left( \begin{array}{cc} rI_n &\quad 0 \\ 0  &\quad 2I_m \end{array}
\right). \]
We know that  $L(G^{00+}) = D(G^{00+}) -A(G^{00+})$ and
$ L(\lambda ,G^{00+}) = \det (\lambda I_n - L(G^{00+}))$.
Therefore
\begin{eqnarray*}
 L(\lambda ,G^{00+})&=& \left| \begin{array}{cc}
 (\lambda  -r)I_n &\quad Q
 \\[1ex]
  Q^{\top}  &\quad (\lambda
-2)I_m
\end{array} \right|.
\end{eqnarray*}
Clearly, it is sufficient to prove our claim for
$\lambda \ne 2$.
Now using Lemmas \ref{lemABCD}
we obtain:
$$L(\lambda ,G^{00+}) =
 |(\lambda - 2)I_m| \times |(\lambda - r)I_n -\frac{1}{\lambda -2} Q I_m Q^{\top}|.$$
By Lemma \ref{lemQ,Qt} $(a1)$, $QQ^{\top} = rI_n + A$.
Therefore
\\[1ex]
\indent
$L( \lambda ,G^{00+})  =  (\lambda  - 2)^{m-n} \times |B|$, where
$B = (\lambda  - 2)( \lambda  - r)I_n - (rI_n + A).$
\\[1ex]
Obviously, $|B|$ equals the product of its  eigenvalues.
By Lemma \ref{Corfg,AJ}, the eigenvalues of $B$ are
\\[0.5ex]
\indent
$\sigma_n = (\lambda  - 2)(( \lambda  - r)) - (r + r) =
\lambda ( \lambda - r -2)$ and
\\[1ex]
\indent
$\sigma_i =
(\lambda  - 2)( \lambda  - r) - (r + r - \lambda _i)
= \lambda ^2 - (r +2) \lambda +  \lambda _i $
for $i = 1,2, \ldots , n-1$.
\\[1ex]
 Since $|B|=\prod_{i=1}^n \sigma_i$,  we have:
 \\[1.5ex]
\indent
 $L( \lambda ,G^{00+}) =
 \lambda (\lambda  - r - 2)(\lambda  - 2)^{m-n}
 \prod_{i=1}^{n-1} \{
\lambda ^2 - (r +2) \lambda +  \lambda _i
   \}.$
 \ep
\begin{corollary}
\label{Sp(00+)}
If $G$ is an $r$-regular graph with $n$ vertices and $m$ edges, then $G^{00+}$ has $m-n$ Laplacian eigenvalues equal to $2$ and the following $2n$  Laplacian eigenvalues
\[ \frac{1}{2}(~r+2\pm \sqrt{(r+2)^2-4\lambda _i}~),\  i=1,2,\cdots,n .\]
\end{corollary}
\begin{corollary}
\label{t(G00+)}
Let $G$ be an $r$-regular graph with $n$ vertices and $m$ edges. Then
\[t(G^{00+}) = \frac{n}{m+n}(r + 2) 2^{m - n} t(G).\]
\end{corollary}
Similarly, we can prove the following theorem.
\begin{Theorem}\label{11+}
Let $G$ be an $r$-regular graph with $n$ vertices and
$m$ edges.
Then
\\[1.5ex]
$(a1)$
$ L( \lambda ,G^{10+}) =
 \lambda (\lambda  - r - 2)(\lambda  - 2)^{m-n}
 \prod_{i=1}^{n-1} \{(\lambda -2)(\lambda - n -r) -2r + \lambda_i \}$,
 \\[1.5ex]
$(a2)$
 $ L( \lambda ,G^{01+}) =
 \lambda (\lambda  - r - 2)(\lambda - m - 2)^{m-n}
 \prod_{i=1}^{n-1} \{(\lambda - r)(\lambda -m -2) -2r
 + \lambda _i \}$ and
\\[1.5ex]
$(a3)$
$L( \lambda ,G^{11+}) =
 \lambda (\lambda  - r - 2)(\lambda - m - 2)^{m-n}
 \prod_{i=1}^{n-1} \{(\lambda - n -r)(\lambda - m -2) -2r + \lambda_i \}.$
\end{Theorem}

\subsection{ Laplacian spectra of $G^{xyz}$ for $z = +$ and $|\{x,y\} \cap \{0,1\}| = 1$}

\begin{Theorem}\label{+0+}
Let $G$ be an $r$-regular graph with $n$ vertices and
$m$ edges.
Then
\[L(\lambda,G^{+0+}) = \lambda (\lambda - r - 2)(\lambda - 2)^{m-n}
\prod_{i=1}^{n-1} \{ (\lambda - 2)(\lambda - r - \lambda_i) - 2r + \lambda_i \}
\]
or equivalently,
\[L(\lambda ,G^{+0+}) = (\lambda - 2)^{m-n}
(\lambda -3)^n
L(\frac{\lambda ^2 - \lambda (r + 2)}{\lambda -3} ,G).\]

\end{Theorem}
\bp
The adjacency matrix and the degree matrix of
$G^{+0+}$ are \[ A(G^{+0+})= \left( \begin{array}{cc} A &\quad
Q\\[1ex]
 Q^{\top}  &\quad 0
\end{array} \right) \hbox{\ and \ } D(G^{+0+})=
\left( \begin{array}{cc} 2rI_n &\quad 0
\\[1ex]
0  &\quad 2I_m
\end{array} \right). \]
Therefore
\begin{eqnarray*}
 L(\lambda ,G^{+0+} )  &=& \left| \begin{array}{cc} (\lambda  -2r)I_n+A &\quad Q
\\[1.3ex] Q^{\top}  &\quad (\lambda -2)I_m
\end{array} \right|.
\end{eqnarray*}
Clearly, it is sufficient to prove our claim for
$\lambda \ne 2$.
Now using Lemmas \ref{lemABCD}
we obtain:
$$ L(\lambda ,G^{+0+} ) =
(\lambda -2)^m   |(\lambda -2r)I_n+A-\frac{1}{\lambda -2}QQ^{\top}|.$$
By Lemma \ref{lemQ,Qt} $(a1)$, $QQ^{\top} = rI_n + A$.
Therefore
\\[1ex]
\indent
$ L(\lambda ,G^{+0+} ) =
(\lambda  - 2)^{m-n}~{B}$, where
$B =
(\lambda - 2)(\lambda -2r)I_n + (\lambda - 2)A -(r I_n + A).$
\\[1ex]
Obviously, $|B|$ equals  the product of its  eigenvalues.
By Lemma \ref{Corfg,AJ}, the eigenvalues of $B$ are
\\[1ex]
\indent
$\sigma_n = (\lambda - 2)(\lambda - 2r) + (\lambda - 2)r -(r  + r) =
\lambda (\lambda  - r - 2)$ and
\\[1ex]
\indent
$\sigma_i =
(\lambda - 2)(\lambda - r -  \lambda _i)  -2r  + \lambda _i$ for $i = 1,2 \ldots , n-1$.
 \\[1.3ex]
Since $|B|=\prod_{i=1}^n \sigma_i$,  we have:
 \\[1.5ex]
 \indent
 $L( \lambda ,G^{+0+}) =
\lambda (\lambda  - r - 2) (\lambda  - 2)^{m-n}
\prod_{i=1}^{n-1} \{ (\lambda - 2) (\lambda - r - \lambda_i) - 2r + \lambda_i \}.$
\ep
\begin{corollary}
\label{Sp(+0+)}
If $G$ is an $r$-regular graph with $n$ vertices and $m$ edges, then $G^{+0+}$ has $m-n$ Laplacian eigenvalues equal to $2$ and the following $2n$  Laplacian eigenvalues
\[ \frac{1}{2}(~r +2 +\lambda _i\pm \sqrt{(r+2+\lambda _i)^2 -12\lambda _i}~),\  i=1,2,\cdots, n .\]
\end{corollary}
\begin{corollary}
\label{t(G+0+)}
Let $G$ be an $r$-regular graph with $n$ vertices and $m$ edges. Then
\[t(G^{+0+}) = \frac{n}{m+n}(r + 2) 2^{m - n} 3^{n-1} t(G),\]
and so
\[t(G^{+0+}) = 3^{n-1} t(G^{00+}).\]
\end{corollary}
\begin{Theorem}
\label{0++}
Let $G$ be an $r$-regular graph with $n$ vertices and $m$ edges.
Then
\[L(\lambda , G^{0++}) = \lambda
(\lambda - r- 2) (\lambda - 2r - 2)^{m-n}
\prod_{i=1}^{n-1}\{ (\lambda - r)(\lambda - 2 - \lambda _i) - 2r + \lambda _i \}
\]
or equivalently,
\[L(\lambda , G^{0++}) =
(\lambda -  r - 1)^n(\lambda - 2r - 2)^{m - n}
L(\frac{\lambda ^2 -  (r + 2)\lambda}{\lambda- r - 1}, G).\]
\end{Theorem}
\bp
The adjacency matrix and the degree matrix of
$G^{0++}$ are \[ A(G^{0++})= \left( \begin{array}{cc} 0 &\quad
Q
\\[1ex]
 Q^{\top}  &\quad A(G^l)
\end{array} \right) \hbox{\ and \ } D(G^{0++})= \left( \begin{array}{cc} rI_n &\quad 0
\\[1ex]
 0  &\quad 2rI_m
\end{array} \right). \]
Since by Lemma \ref{lemQ,Qt} $(a2)$,
$A(G^l) =  Q^{\top}Q -  2I_m$, we have:
\begin{eqnarray*}
 L(\lambda , G^{0++})  &=&
 \left| \begin{array}{cc} (\lambda - r)I_n &~~ \quad Q
 \\[1.5ex]
 Q^{\top}  &~~ \quad (\lambda - 2r - 2)I_m + Q^{\top}Q
\end{array} \right|
\\[3ex]
 & = & \left| \begin{array}{cc}
 (\lambda - r)I_n &~~ \quad Q
 \\[1.5ex]
 (r + 1 - \lambda )Q^{\top}
  &~~ \quad (\lambda  - 2r - 2)I_m
\end{array} \right|.
\end{eqnarray*}
Clearly, it is sufficient to prove our claim for
$\lambda \ne 2r + 2$.
Using Lemmas \ref{lemABCD}
we obtain:
\begin{eqnarray*}
 L(\lambda , G^{0++})  &=& (\lambda - 2r - 2)^m|(\lambda -r)I_n-\frac{r+1-\lambda }{\lambda -2r-2}QQ^{\top}|.
\end{eqnarray*}
By Lemma \ref{lemQ,Qt} $(a1)$,
$QQ^{\top} = rI_n + A$.
Therefore
\\[1ex]
\indent
$ L(\lambda , G^{0++})  = (\lambda - 2r - 2)^{m-n}~|B|$, where
$B = ( \lambda - r)( \lambda - 2r - 2)I_n - (r + 1- \lambda )(r I_n + A)$.
\\[1ex]
Obviously, $|B|$ equals  the product of its  eigenvalues.
By Lemma \ref{Corfg,AJ}, the eigenvalues of $B$ are
\\[0.7ex]
\indent
$\sigma_n = (\lambda - r)( \lambda - 2r - 2) - (r + 1- \lambda )(r  + r) =
 \lambda (\lambda - r - 2)$ and
 \\[1ex]
\indent
$\sigma_i =
(\lambda - r)(\lambda - 2 -\lambda_i) - 2r + \lambda_i$ for $i = 1, 2, \ldots , n-1$.
 \\[1.3ex]
Since $|B|=\prod_{i=1}^n \sigma_i$,  we have:
 \\[1.5ex]
 \indent
$ L(\lambda , G^{0++}) =
\lambda (\lambda - r - 2) (\lambda - 2r -2)^{m-n}
\prod_{i=1}^{n-1} \{(\lambda - r)(\lambda - 2 -\lambda_i) - 2r + \lambda_i\}$.
\ep
\begin{corollary}
\label{Sp(0++)}
If $G$ is an $r$-regular graph with $n$ vertices and $m$ edges, then $G^{0++}$ has $m-n$ Laplacian eigenvalues equal to $2r+2$ and the following $2n$  Laplacian eigenvalues
\[ \frac{1}{2}(~r + 2 + \lambda _i\pm \sqrt{
\lambda _i^2 - 2r \lambda _i + r^2 + 4r + 4}~),~~i = 1, 2, \cdots , n .\]
\end{corollary}

\begin{corollary}
\label{t(G0++)}
Let $G$ be an $r$-regular graph with $n$ vertices and $m$ edges. Then
\[t (G^{0++}) =
\frac{n}{m+n} 2^{m - n}(r + 1)^{m-1}(r+2)   t(G).\]
\end{corollary}

Theorems \ref{00+}  and \ref{0++} and Corollaries
\ref{t(G00+)}, \ref{t(G+0+)}, and \ref{t(G0++)}
coincide with the corresponding Theorems \ref{G0++}, \ref{G00+},  and \ref{trees} by Kelmans.
\\[1.5ex]
\indent
The remaining situations when $z = +$ and
$|\{x,y\} \cap \{0,1\}| = 1$ can be considered similarly.
Here is the  list of the Laplacian characteristic polynomials of  $G^{xyz}$ for these situations.
\begin{Theorem}
\label{otherxy01z+}
Let $G$ be an $r$-regular graph with $n$ vertices and $m$ edges. Then
\\[1ex]
$ L( \lambda ,G^{0-+}) =
 \lambda (\lambda  - r - 2)(\lambda - m + 2r - 2)^{m-n}
 \prod_{i=1}^{n-1} \{(\lambda -r)(\lambda - m - 2 - \lambda_i) -2r + \lambda_i \},  $
\\[1ex]
$ L( \lambda ,G^{-0+}) =
 \lambda (\lambda  - r - 2)(\lambda  - 2)^{m-n}
 \prod_{i=1}^{n-1} \{(\lambda - n -r + \lambda_i)(\lambda -2) -2r + \lambda_i \},  $
\\[1ex]
$ L( \lambda ,G^{1-+}) =
 \lambda (\lambda  - r - 2)(\lambda -m + 2r - 2)^{m-n}
 \prod_{i=1}^{n-1} \{(\lambda - n -r)(\lambda - m -2 + \lambda_i) -2r + \lambda_i \},  $
 \\[1ex]
$  L( \lambda ,G^{-1+}) =
 \lambda (\lambda  - r - 2)(\lambda - m - 2)^{m-n}
 \prod_{i=1}^{n-1} \{(\lambda - n -r + \lambda_i)(\lambda - m -2) -2r + \lambda_i \}, $
\\[1ex]
$ L( \lambda ,G^{+1+}) =
 \lambda (\lambda  - r - 2)(\lambda - m - 2)^{m-n}
 \prod_{i=1}^{n-1} \{(\lambda - r - \lambda_i)(\lambda - m -2) -2r + \lambda_i \}   $ and
\\[1ex]
$ L( \lambda ,G^{1++}) =
 \lambda (\lambda  - r - 2)(\lambda  - 2r - 2)^{m-n}
 \prod_{i=1}^{n-1} \{(\lambda - n -r)(\lambda -2 - \lambda_i) -2r + \lambda_i \}.  $
\end{Theorem}

\subsection{ Laplacian spectra of $G^{xyz}$ for $z \in\{1, -\} $}

\indent

Using the Reciprocity Theorem \ref{reciprocity}
 it is easy to find from Theorems \ref{xy0}, \ref{00+},
 \ref{11+},
 \ref{+0+}, \ref{0++} and \ref{otherxy01z+}
the Laplacian characteristic polynomial of every $G^{xyz}$ with $\{x, y\} \cap \{0, 1\} \neq \emptyset$ and $z \in \{ 1, -\}$ (see Appendix).

\section{Laplacian spectra of $G^{xyz}$ with
$x,y,z \in \{+, -\}$}
\label{x,y,z-in-(+,-)}

\indent

 In this section we will describe
the Laplacian characteristic polynomials and the Laplacian spectra  of transformations $G^{xyz}$ of an $r$-regular graph $G$ for
$x,y,z \in \{+, -\}$ in terms of the Laplacian spectrum of
$G$, $r$, $v(G) = n$, $r$ (and $e(G) = m = \frac{1}{2}r n$).
\begin{Theorem}
\label{th+++}
Let $G$ be an $r$-regular graph with $n$ vertices and $m$ edges.
Then
\[L(\lambda ,G^{+++} ) =
\lambda (\lambda  - r - 2)(\lambda - 2r -2)^{m - n}
\prod_{i=1}^{n - 1}
\{ (\lambda - r - \lambda_i)(\lambda - 2 - \lambda_i) - 2r + \lambda_i \}. \]
\end{Theorem}
\bp
The adjacency matrix and the degree matrix of
$G^{+++}$ are
\[A(G^{+++})= \left( \begin{array}{cc} A &\quad Q
\\[1ex]
 Q^{\top}  &\quad
A(G^l)
\end{array}
\right)~~ \hbox{\ and\ }~~ D(G^{+++})= \left( \begin{array}{cc} 2rI_n &\quad 0
\\[1ex]
 0  &\quad
2rI_m
\end{array} \right).\]
Since $L(G^{+++}) = D(G^{+++}) - A(G^{+++})$,
$L( \lambda ,G^{+++}) = \det (\lambda I_n - L(G^{+++}))$, and by Lemma \ref{lemQ,Qt} $(a2)$,
$A(G^l) = Q^{\top}Q - 2I_m$,
we have:
\begin{eqnarray*}
 L( \lambda ,G^{+++}) &=& \left| \begin{array}{cc} ( \lambda  - 2r) I_n + A &~~ \quad Q
 \\[2.5ex]
  Q^{\top}  &~~ \quad ( \lambda  - 2r -2) I_m + Q^{\top}Q
\end{array} \right|
\\[4ex]
& = & \left| \begin{array}{cc} ( \lambda  - 2r)I_n + A &~~ \quad Q
\\[2.5ex]
 - Q^{\top}(( \lambda  - 2r - 1)I_n + A)  &~~ \quad
( \lambda  - 2r - 2)I_m
\end{array} \right|.
\end{eqnarray*}
Clearly, it is sufficient to prove our claim for
$\lambda \ne 2r +2$.
Using Lemmas \ref{lemABCD}
we obtain:
$$L( \lambda ,G^{+++})  =
|( \lambda - 2r - 2)I_m|\times |( \lambda -2r)I_n + A +
(\lambda - 2r - 2)^{-1}QQ^{\top}(( \lambda -2r-1)I_n + A)|.
$$
By Lemma \ref{lemQ,Qt} $(a1)$, $QQ^{\top} = rI_n + A$. Therefore
$L( \lambda ,G^{+++})  =  (\lambda - 2r - 2)^{m-n} \times |B|$, where
$$B = |( \lambda - 2r - 2)(( \lambda - 2r) I_n + A) +
(rI_n + A) (( \lambda - 2r - 1)I_n + A)|.$$
Obviously,  $|B|$ equals  the product of its eigenvalues.
By Lemma \ref{Corfg,AJ},
the eigenvalues of $B$ are
\\[1ex]
$\sigma_i =
( \lambda - 2r - 2)(( \lambda - 2r) + r - \lambda _i) +
(r + r - \lambda _i) (( \lambda - 2r - 1) + r - \lambda _i)
\\[0.5ex]
~~~~=(\lambda -  r - \lambda _i) ( \lambda  - 2  - \lambda _i)
- 2 r + \lambda _i$
for $i = 1, 2, \ldots , n$.
\\[1ex]
Since $|B|=\prod_{i=1}^n \sigma_i$
and $\lambda _n = 0$,
we have:
\\[1.5ex]
\indent
$L(\lambda ,G^{+++} ) =
\lambda (\lambda  - r - 2)(\lambda - 2r -2)^{m - n}
\prod_{i=1}^{n - 1}
\{ (\lambda - r - \lambda_i)(\lambda - 2 - \lambda_i) - 2r + \lambda_i \}. $
\ep

\begin{Theorem}
\label{th+-+}
Let $G$ be an $r$-regular graph with $n$ vertices and $m$ edges.
Then
\[L(\lambda ,G^{+-+}) = \lambda (\lambda - r - 2)
(\lambda - m + 2r - 2)^{m - n}
\prod_{i=1}^{n-1}\{ (\lambda - r - \lambda_i) (\lambda - m - 2 + \lambda _i) - 2r + \lambda_i\}.
\]
\end{Theorem}

\bp
From the definition of $G^{+-+}$, we have:
\[A(G^{+-+}) = \left(
                 \begin{array}{cc}
                   A &~~ Q \\[1.5ex]
                   Q^\top &~~ J_{m m} - I_m - A(G^l) \\
                 \end{array}
               \right).
\]
For every $z \in V(G^ {+-+})$, $d(z,G^{+-+}) = 2r$ if $z\in V(G)$ and $d(z,G^{+-+}) = 2+m-1-(2r-2) = m - 2r + 3$ if
$z \in E(G)$.
Therefore,
\[ D(G^{+-+}) = \left(
                                                              \begin{array}{cc}
                                                                2rI_n &~~ 0
                                                                \\[1.5ex]
                                                                0 &~~ (m - 2r + 3)I_m \\
                                                              \end{array}
                                                            \right).
\]
Then $L(\lambda,G^{+-+}) = |M|$, where
\begin{eqnarray*}
 M &=& \left(
                 \begin{array}{cc}
                   (\lambda - 2r)I_n + A &~~~ Q \\[2ex]
                   Q^\top &~~~~ (\lambda - m + 2r - 2) I_m + J_{mm} - Q^\top Q \\
                 \end{array}
               \right).
\end{eqnarray*}
\\
By Lemma \ref{lemJ} $(a4)$,  $J_{mn}Q = 2J_{mm}$.
Hence multiplying the first row of the block matrix $M$ by
 $Q^\top - \frac{1}{2} J_{mn}$ and adding the result to the second row of $M$, we obtain a new matrix
\begin{eqnarray*}
 M' &=& \left(
                 \begin{array}{cc}
                   (\lambda - 2r)I_n + A &~~~~ Q \\[2ex]
 Q^\top ((\lambda - 2r+ 1)I_n + A) - \frac{1}{2}J_{mn}((\lambda - 2r)I_n + A) &~~~~ (\lambda - m + 2r - 2) I_m\\
                 \end{array}
               \right).
\end{eqnarray*}
\\
Clearly, $L(\lambda,G^{+-+}) = |M| = |M'|$.
Obviously, it is sufficient to prove our claim for
$\lambda \ne m -2r +2$.
Now using Lemmas \ref{lemABCD}
we obtain:
\\[1ex]
\indent
$L(\lambda,G^{+-+}) = (\lambda - m + 2r - 2)^{m-n} \times |B|$,
where
\\[1ex]
$B = (\lambda - m + 2r - 2)((\lambda - 2r)I_n + A)
- QQ^\top((\lambda - 2r + 1)I_n + A) +\frac{1}{2}QJ_{mn}((\lambda - 2r)I_n + A)$.
\\[1ex]
By Lemma \ref{lemQ,Qt} $(a1)$, $QQ^\top = rI_n + A$ and
by Lemma \ref{lemJ} $(a2)$, $QJ_{mn} = r J_{nn}$. Therefore
\\[1ex]
$B = (\lambda - m + 2r - 2)((\lambda - 2r)I_n + A)
- (rI_n + A)((\lambda - 2r + 1)I_n + A) +\frac{r}{2}
J_{nn}((\lambda - 2r) + A).$
\\[1ex]
Obviously,  $|B|$ equals the product of its eigenvalues and $nr = 2m$.
By Lemma \ref{Corfg,AJ},
the eigenvalues of $B$ are
\\[1ex]
$\sigma_n = (\lambda - m + 2r - 2)(\lambda - 2r + r) -2r(\lambda -2r + 1 + r) + \frac{nr}{2}(\lambda - 2r +r)
= \lambda^2 - \lambda(r+2)$ and
\\[1ex]
$\sigma_i =  (\lambda - m + 2r - 2)(\lambda -2r + r -\lambda_i)- (2r -\lambda_i)(\lambda  -2r + r + 1 -\lambda_i)
\\[0.5ex]
~~~~= (\lambda - r - \lambda_i) (\lambda - m - 2 + \lambda _i) - 2r + \lambda_i$
for $i = 1, 2 ,\cdots, n - 1.$
\\[1ex]
 Since $|B|=\prod_{i=1}^n \sigma_i$,  we have:
 \\[1.5ex]
$L(\lambda ,G^{+-+}) = \lambda (\lambda - r - 2)
(\lambda - m + 2r - 2)^{m - n}
\prod_{i=1}^{n-1}\{ (\lambda - r - \lambda_i) (\lambda - m - 2 + \lambda _i) - 2r + \lambda_i\}.$

\ep
\begin{Theorem}
\label{th-++}
Let $G$ be an $r$-regular graph with $n$ vertices and $m$ edges.
Then
\[ L(\lambda, G^{-++}) = \lambda(\lambda-r-2)
(\lambda - 2r - 2)^{m-n}
\prod_{i=1}^{n-1} \{(\lambda -n -r +\lambda_i)(\lambda  - 2  - \lambda_i)- 2r + \lambda_i \}. \]
\end{Theorem}

\bp
From the definition of
$G^{-++}$ we have:
\[ A(G^{-++})= \left( \begin{array}{cc} J_{nn} - I_n -A &\quad Q \\
Q^{\top}  & \quad A(G^l)
\end{array} \right). \]
For every vertex $z$ in $G^{-++}$, we have: $d(z,G^{-++}) = n-1$ if $z\in V(G)$ and $d(z,G^{-++}) = 2r$ if $z\in E(G)$, and so  \[D(G^{-++}) = \left( \begin{array}{cc} (n-1)I_n &\quad 0 \\
0  & \quad 2rI_m
\end{array} \right). \]
By  Lemma \ref{lemQ,Qt} $(a2)$,
$A(G^l) = 2I_m - Q^{\top}Q $.
Therefore

\begin{eqnarray*}
  L(\lambda , G^{-++}) &=& \left|
                             \begin{array}{cc}
                               (\lambda -n +1)I_n + J_{nn} - I_n -A & Q \\[1.5ex]
                               Q^{\top} &~~ (\lambda- 2r)I_m + Q^{\top}Q -2I_m \\
                             \end{array}
                           \right| .
\end{eqnarray*}
Multiplying the first row of the above block matrix by
$- Q^{\top}$ and adding the result to the second row we obtain:
\begin{eqnarray*}
  L(\lambda , G^{-++})  &=&  \left|
                       \begin{array}{cc}
                         (\lambda - n)I_n + J_{nn} -A &~~ Q \\[1.5ex]
                         Q^{\top} - Q^{\top}((\lambda - n)I_n + J_{nn} -A) &~~ (\lambda - 2r -2)I_m \\
                       \end{array}
                     \right| .
\end{eqnarray*}

Clearly, it is sufficient to prove our claim for
$\lambda \ne 2r + 2$.
Using Lemmas \ref{lemABCD}
we obtain:
\[
\begin{array}{rl}
 L(\lambda,G^{-++}) =& (\lambda - 2r - 2)^{m-n}\times |B|,
  \end{array} \]
 where
 \[
\begin{array}{rl}
B = & (\lambda - 2r - 2)((\lambda - n)I_n + J_{nn} - A) -
 QQ^\top((1-\lambda + n)I_n - J_{nn} + A).
 \end{array} \]
 Since by Lemma \ref{lemQ,Qt} $(a1)$, $QQ^\top = D(G) + A(G) = rI_n + A$, we have:
  \[
\begin{array}{rl}
B = & (\lambda - 2r - 2)((\lambda - n)I_n + J_{nn} - A )-
(rI_n + A)((1-\lambda + n)I_n - J_{nn} + A).
 \end{array} \]
 Obviously, $|B|$ is the product of the eigenvalues of $B$.
 By Lemma \ref{Corfg,AJ},
the eigenvalues of $B$ are
\\[1ex]
$\sigma_n = (\lambda - 2r - 2)(\lambda - n + n  - r )-
(r + r)(1-\lambda + n - n + r)  =
\lambda (\lambda - r - 2)$
and
\\[1ex]
$\sigma_i =
(\lambda - 2r - 2) (\lambda - n   - r + \lambda _i ) -
(r + r - \lambda _i)((1 - \lambda + n) + r - \lambda _i)
\\[0.5ex]
~~~~=(\lambda - 2 - \lambda _i )
( \lambda - n - r + \lambda _i ) - 2r + \lambda _i $
for $i = 1,2, \ldots , n-1$.
\\[1ex]
 Since $|B|=\prod_{i=1}^n \sigma_i$,  we have:
\\[1.5ex]
$L(\lambda,G^{-++}) = \lambda(\lambda-r-2)
(\lambda - 2r - 2)^{m - n}\prod_{i=1}^{n-1} \{ (\lambda - n - r + \lambda_i)(\lambda  - 2  - \lambda_i)- 2r + \lambda_i \}.$

\ep

Similarly we can prove the following theorem  for
$G^{--+}$:
\begin{Theorem}\label{th--+}
Let $G$ be an $r$-regular graph with $n$ vertices and $m$ edges. Then
$$L(\lambda, G^{--+}) = \lambda ( \lambda - r - 2) (\lambda - m + 2r - 2)^{m-n} \prod_{i=1}^{n-1} \{ (\lambda - n - r + \lambda_i)(\lambda - m  - 2 + \lambda_i) - 2r + \lambda_i \}.$$
\end{Theorem}
Now we can use Reciprocity Theorem \ref{reciprocity} to obtain from
Theorems \ref{th+++} - \ref{th--+} the Laplacian characteristic polynomials of the corresponding complement graphs $G^{xyz}$.

\begin{Theorem}
\label{complements}
Let $G$ be an $r$-regular graph with $n$ vertices and $m$ edges and let $s = n + m$.
Then
\\[1.5ex]
$(a1)$
$
L( \lambda ,G^{---})
=  \lambda (\lambda  - s + r + 2)
(\lambda - s + 2r + 2)^{m - n} \times
\\[1ex]
\indent
~~~~~~~~~~~~~~~~~~~\prod_{i=1}^{n-1}\{(\lambda - s + r + \lambda _i)
(\lambda - s + 2 + \lambda _i) - 2r + \lambda _i\},
$
\\[2.5ex]
$(a2)$
$
L(\lambda , G^{-+-}) =
\lambda (\lambda - s + r + 2)
(\lambda - n - 2r + 2)^{m - n} \times
\\[1ex]
\indent
~~~~~~~~~~~~~~~~~~~\prod_{i=1}^{n-1}
 \{ (\lambda - s  + r + \lambda_i)(\lambda - n + 2 - \lambda_i) - 2r + \lambda_i\},
 $
\\[2.5ex]
$(a3)$
$
L(\lambda , G^{+--})
=  \lambda (\lambda - s + r + 2)(\lambda - s + 2r + 2)^{m-n} \times
\\[1ex]
\indent
~~~~~~~~~~~~~~~~~~~\prod_{i=1}^{n-1}
\{ (\lambda -m + r - \lambda_i)(\lambda - s + 2 + \lambda_i) - 2r + \lambda_i \},
$  and
\\[2.5ex]
$(a4)$
$
L(\lambda ,G^{++-})=
\lambda (\lambda -s +r +2) (\lambda - n -2r + 2)^{m-n}
\times
\\[1ex]
\indent
~~~~~~~~~~~~~~~~~~~\prod_{i=1}^{n-1} \{ (\lambda - m + r - \lambda_i) (\lambda - n + 2 - \lambda_i) - 2r + \lambda_i \}.
$
\end{Theorem}

\bp
As we have mentioned above,
the claims $(a1)$ - $(a4)$ can be easily proven from  Theorems \ref{th+++} - \ref{th--+}, respectively, using   Reciprocity Theorem \ref{reciprocity}.
We give below the proof of claim $(a1)$.
The proofs of the remaining claims $(a2)$ - $(a4)$ are similar.
\\[1ex]
\indent
Since
$G^{+++}$ and $G^{---}$ are complement, we can apply Reciprocity Theorem to obtain
 from Theorem \ref{th+++}:
\\[2ex]
$L( \lambda ,G^{---})  =  (-1)^{s-1} \frac{\lambda }{s - \lambda }
L (s - \lambda , G^{+++})$
\\[2ex]
$ =  (-1)^{s-1} \frac{\lambda }{s - \lambda }
(s - \lambda )(s - \lambda  - r - 2)
(s - \lambda - 2r - 2)^{m - n}
\prod_{i=1}^{n-1}\{(s - \lambda - r - \lambda _i)
(s - \lambda - 2 \lambda _i) - 2r + \lambda _i\}$
\\[2.5ex]
$ =  \lambda (\lambda  - s + r + 2)
(\lambda - s + 2r + 2)^{m - n}
\prod_{i=1}^{n-1}\{(\lambda - s + r + \lambda _i)
(\lambda - s + 2 + \lambda _i) - 2r + \lambda _i\}.
$
\ep
\\[3ex]
\indent
From the above results it follows that
the transformations $G^{xyz}$ have the following common Laplacian spectrum properties.
\begin{Theorem}
\label{z+}
Let $G$ be an $r$-regular graph with $n$ vertices and $m$ edges and
$F = G^{xyz}$, where $z \in \{+, -\}$.
Then $F$ and $F^c$ have, respectively, the Laplacian eigenvalue
\\[1ex]
$(a1)$
$r + 2$ and $ m+ n - r -2$ of multiplicity one if $z = +$,
\\[1ex]
$(a2)$
$2r + 2$ and $m+ n -2r - 2$ of multiplicity $m - n$
 if $(y,z) = (+,+)$,
\\[1ex]
$(a3)$
$m - 2 r + 2$ and $n + 2r - 2$ of multiplicity $m - n$
if $(y,z) =(-,+)$,
\\[1ex]
$(a4)$
$2$ and $m+ n - 2$ of  multiplicity $m-n$ if
$(y,z) = (0, +)$, and
\\[1ex]
$(a5)$
$m + 2$ and $n - 2$ of  multiplicity $m-n$ if $(y,z)= (1,+)$.
\end{Theorem}

The proofs of Theorems
\ref{11+}$(a1)$,
\ref{11+}$(a2)$, and  \ref{th--+} can be found in \cite{DKMarxiv}.

\section{Transformation graphs of  cycles}
\label{examples}

\indent

In this section we first describe some $xyz$-transformations of the 4-cycle and the 5-cycle  and the Laplacian spectra of these transformations.
After that we consider $xyz$-transformations of any cycle and show that some different $xyz$-transformations of the same cycle may be isomorphic.

\begin{figure}[ht]
\begin{center}
\scalebox{0.35}[.35]{\includegraphics{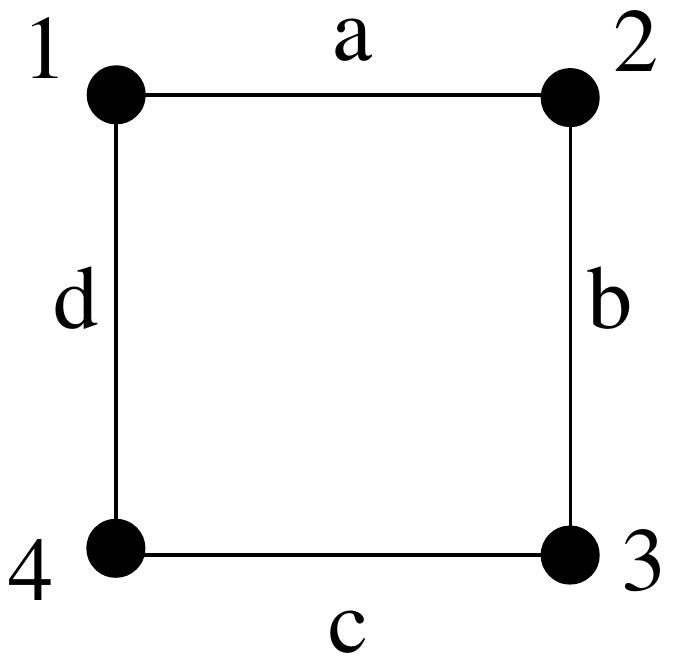}} \ \
\scalebox{0.4}[.4]{\includegraphics{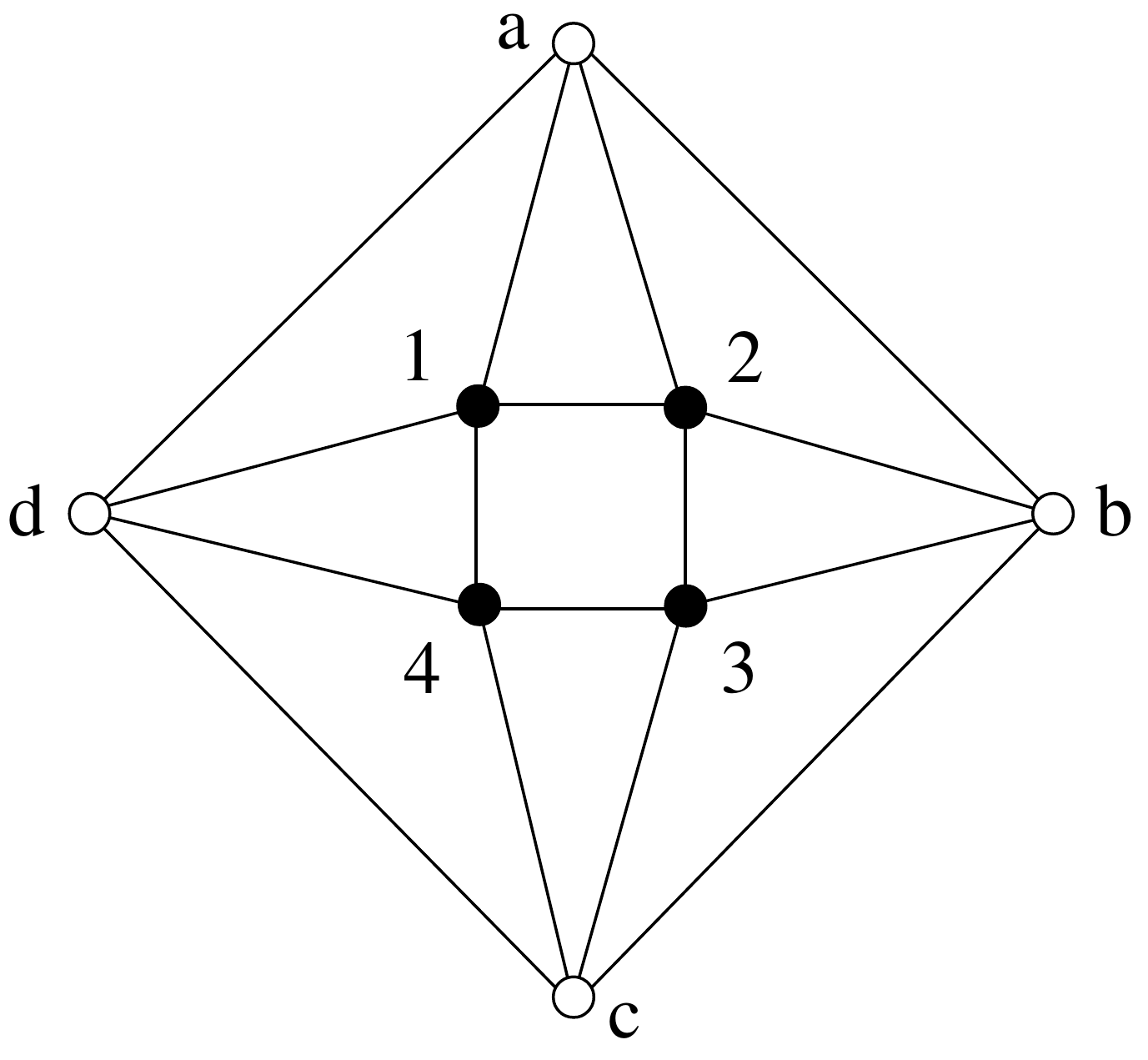}} \ \
\scalebox{0.4}[.4]{\includegraphics{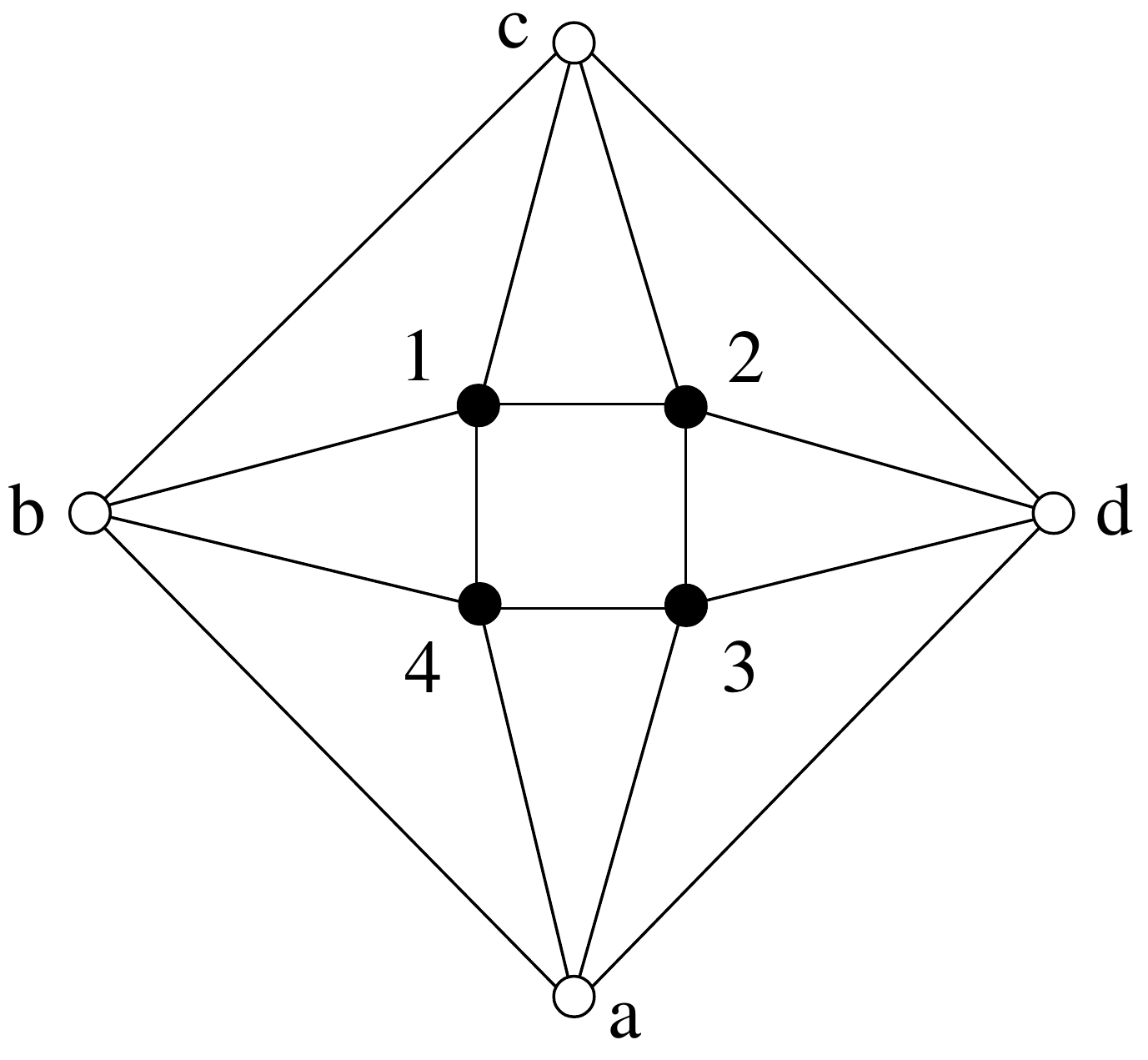}}
\end{center}
\caption{The 4-cycle $C_4$ and its transformations $C_4^{+++}$ (middle) and $C_4^{++-}$ (right)}\label{C4+++}
\end{figure}

\begin{figure}[ht]
\begin{center}
\scalebox{0.4}[.4]{\includegraphics{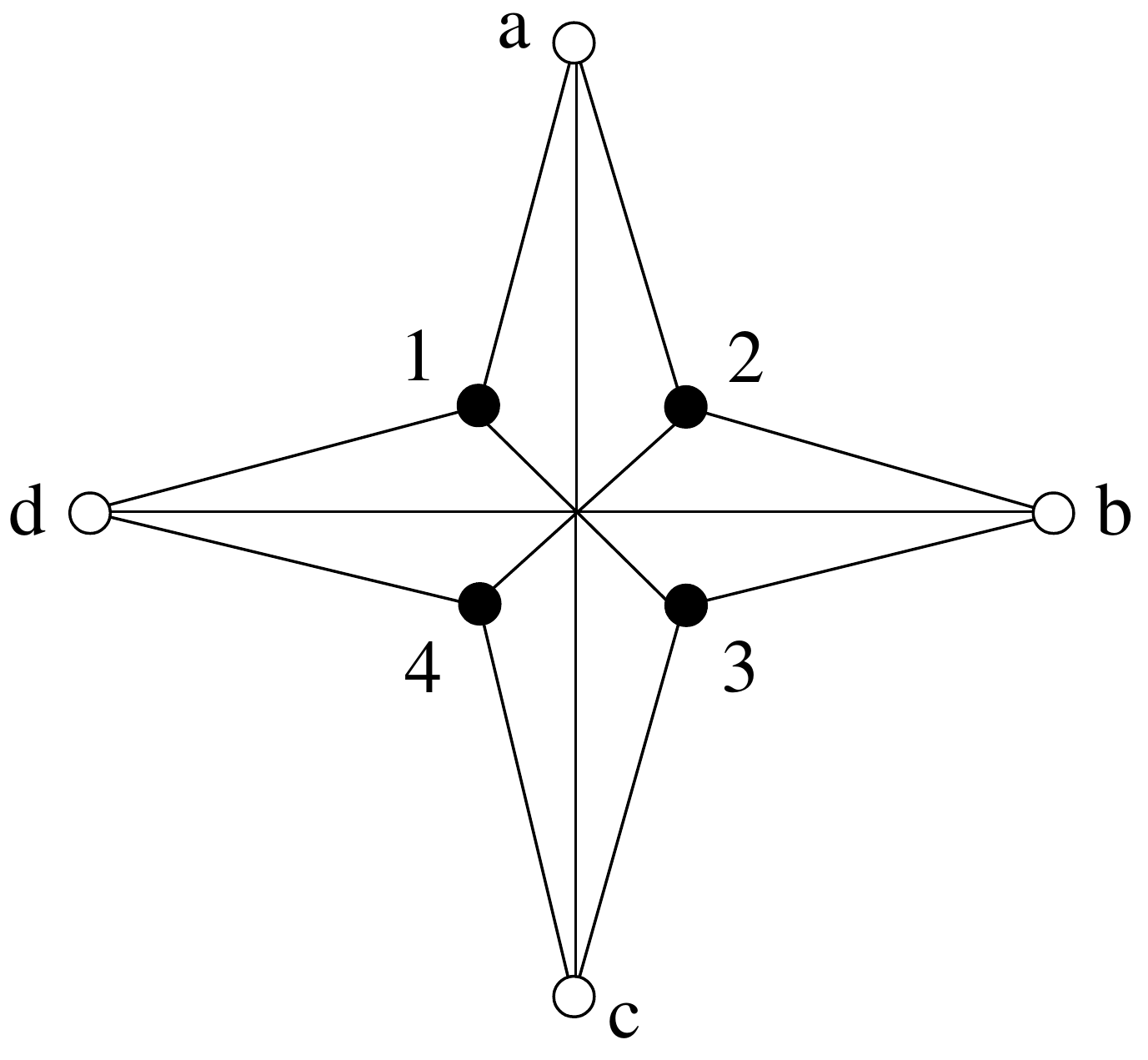}} \quad \quad
\scalebox{0.4}[.4]{\includegraphics{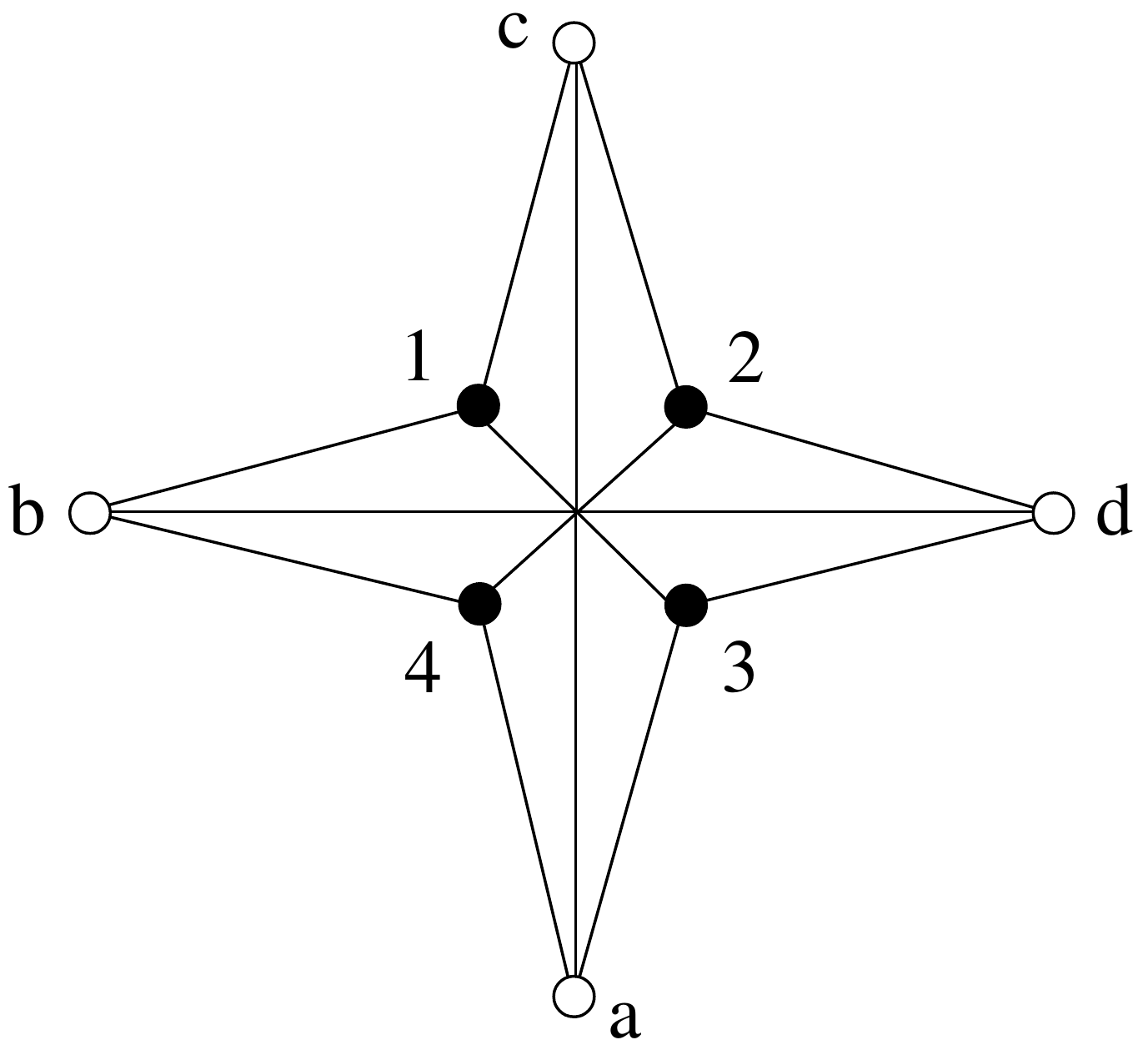}}
\end{center}
\caption{ $C_4^{--+}$ (left) and $C_4^{---}$ (right)}\label{C4--+}
\end{figure}

\begin{figure}[ht]
\begin{center}
\scalebox{0.4}[.4]{\includegraphics{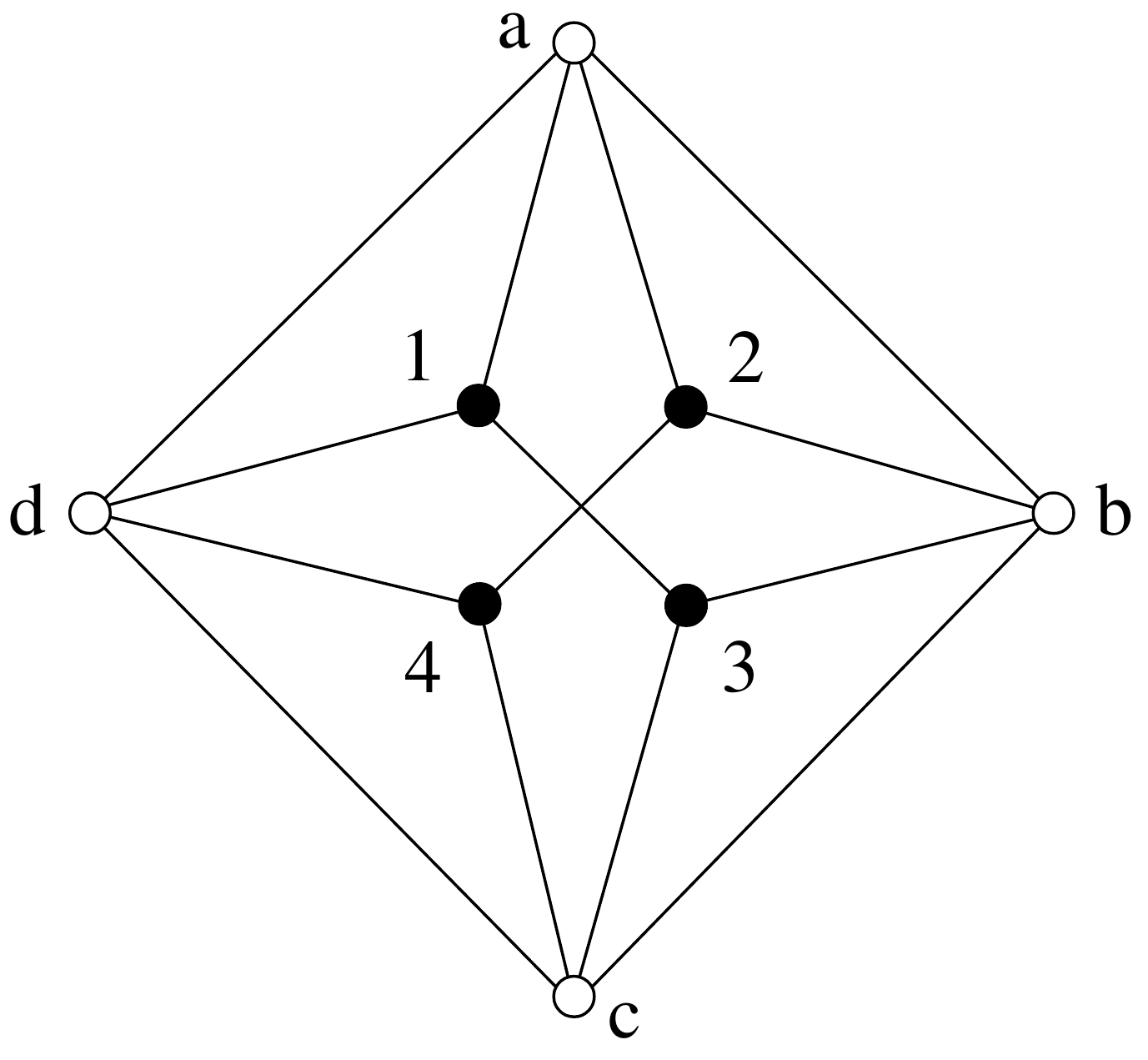}} \quad \quad
\scalebox{0.4}[.4]{\includegraphics{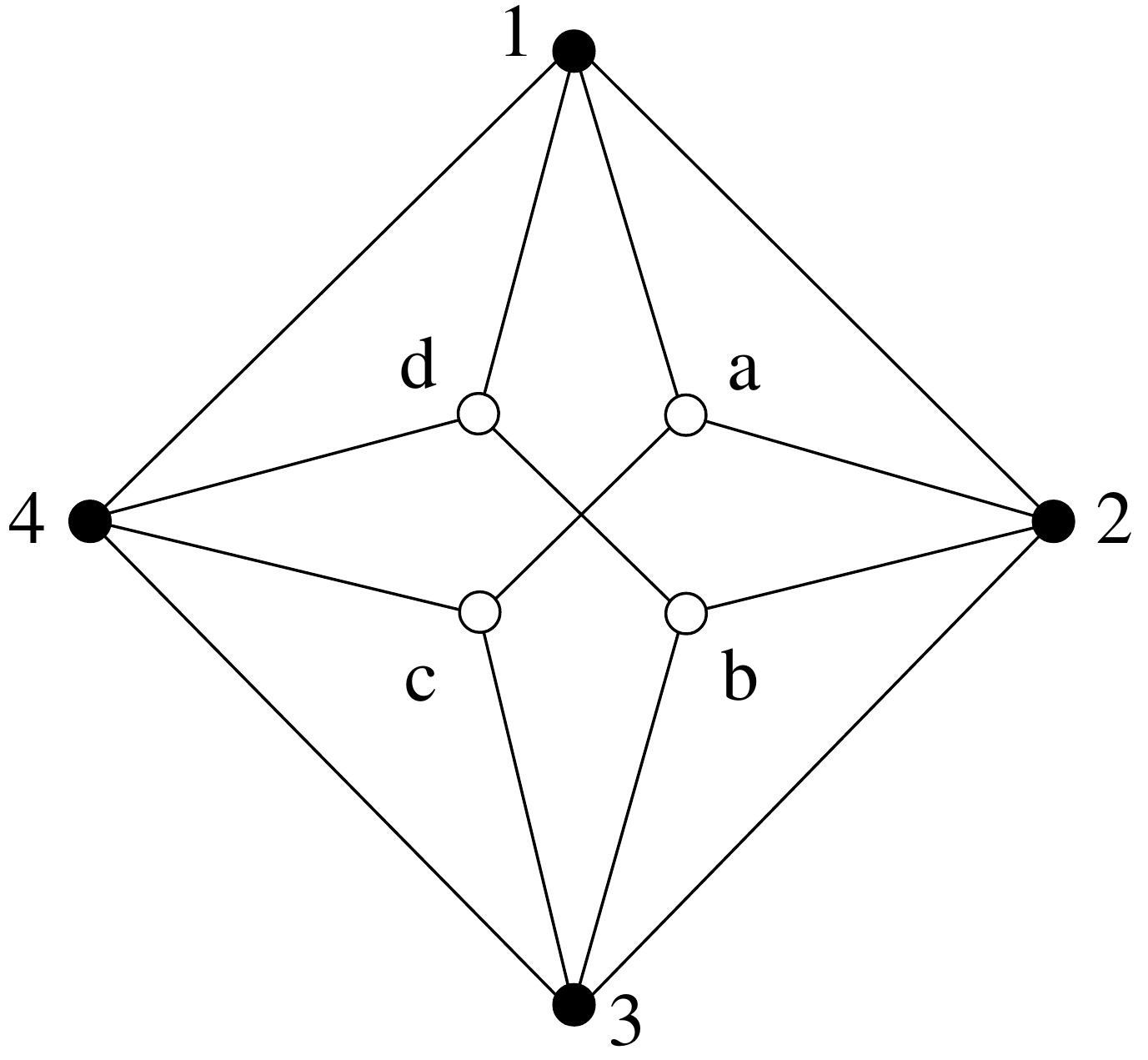}}
\end{center}
\caption{$C_4^{-++}$ (left) and $C_4^{+-+}$ (right)}\label{C4-++/+-+}
\end{figure}
Let $C_n$ be the cycle with $n$ vertices.
It is known (see, for example, \cite{CDS}) that
\\[1ex]
\indent~~~~~~~~~~~~~~~~$Sp(C_n) =
\{2 - 2 ~cos~\frac{2\pi}{n} i: i = 1,\ldots , n\} $.
\\[1ex]
\indent
Let $G$ be a 4-cycle $C_4$. It is easy to see
(and it  follows from the above formula for $C_n$) that
$Sp(C_4) = \{4, 2^{(2)}, 0\}$. Then
\\[1.3ex]
${\bf (c1)}$
$C_4^{+++}$ and $C_4^{++-}$ are isomorphic and
by Theorem \ref{th+++} or \ref{complements}
\\[1ex]
\indent
$Sp(C_4^{+++}) = Sp(C_4^{++-}) = \{6^{(2)}, [4 + \sqrt{2}]^{(2)}, 4, [4 - \sqrt{2}]^{(2)}, 0 \}$,
\\[1.3ex]
${\bf (c2)}$
$C_4^{--+}$ and $C_4^{---}$ are isomorphic and
by Theorem \ref{th--+} or \ref{complements},
\\[1ex]
\indent
$Sp(C_4^{--+} )= = Sp(C_4^{---}) = \{ [4 + \sqrt{2}]^{(2)}, 4, [4 - \sqrt{2}]^{(2)}, 2^{(2)}, 0 \}$,
and
\\[1.3ex]
${\bf (c3)}$
 $C_4^{-++}$ and $C_4^{+-+}$ are isomorphic and by
Theorem \ref{th+-+} or \ref{th-++},
\\[1.5ex]
\indent
$Sp(C_4^{-++}) = Sp(C_4^{+-+}) =
\{6,[4 + \sqrt{2}]^{(2)}, 4, [4 - \sqrt{2}]^{(2)}, 2, 0 \}$.
\begin{figure}[ht]
\begin{center}
\scalebox{0.4}[.4]{\includegraphics{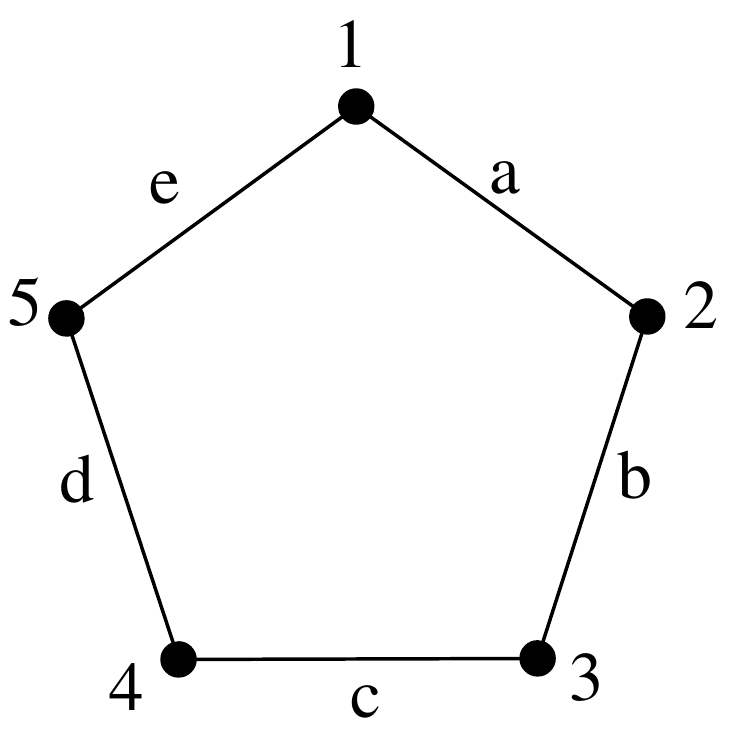}} \ \
\scalebox{0.4}[.4]{\includegraphics{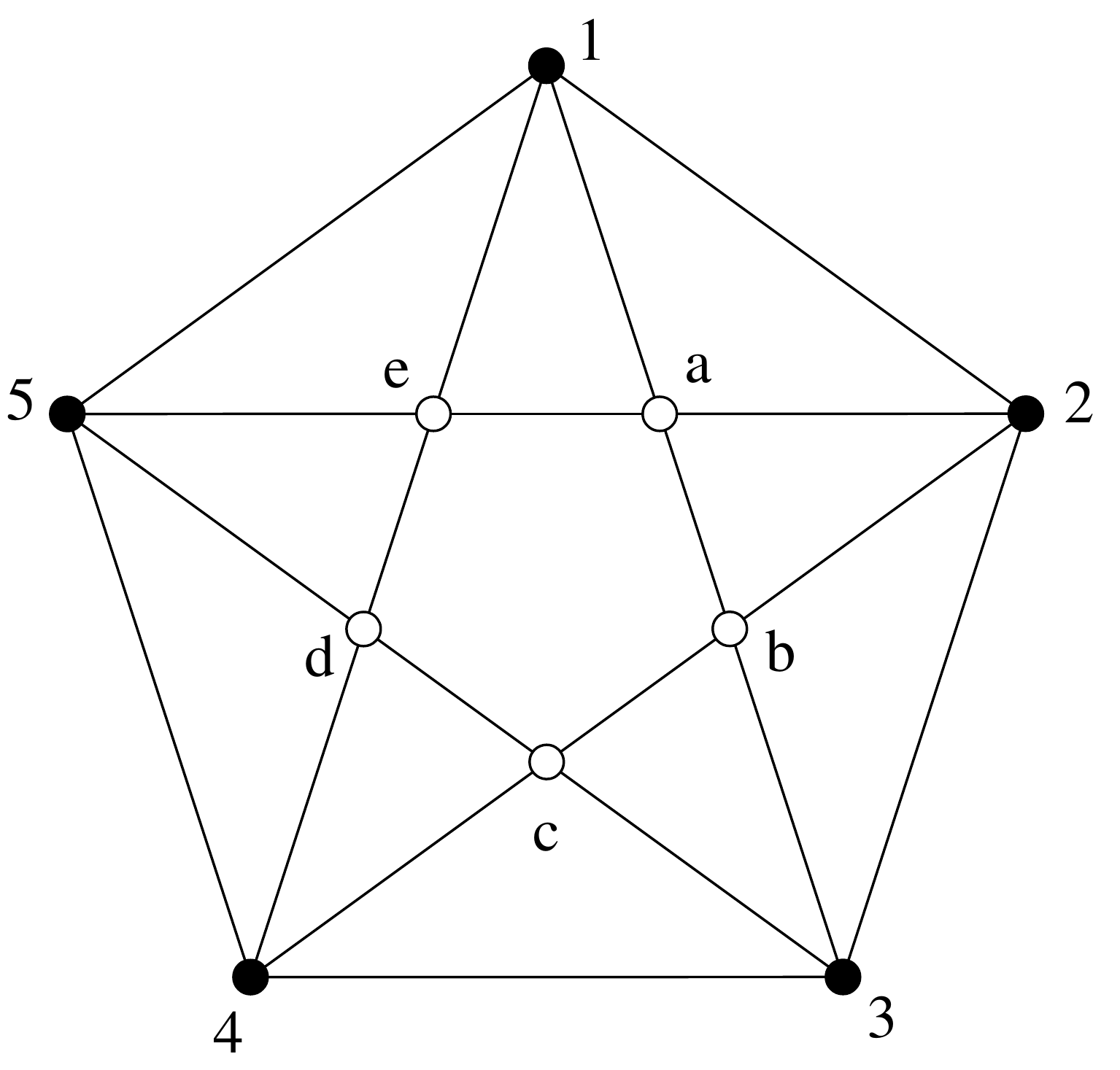}}
\end{center}
\caption{The 5-cycle $C_5$ and its transformation $C_5^{+++}$}\label{C+++}
\end{figure}
\begin{figure}[ht]
\begin{center}
\scalebox{0.4}[.4]{\includegraphics{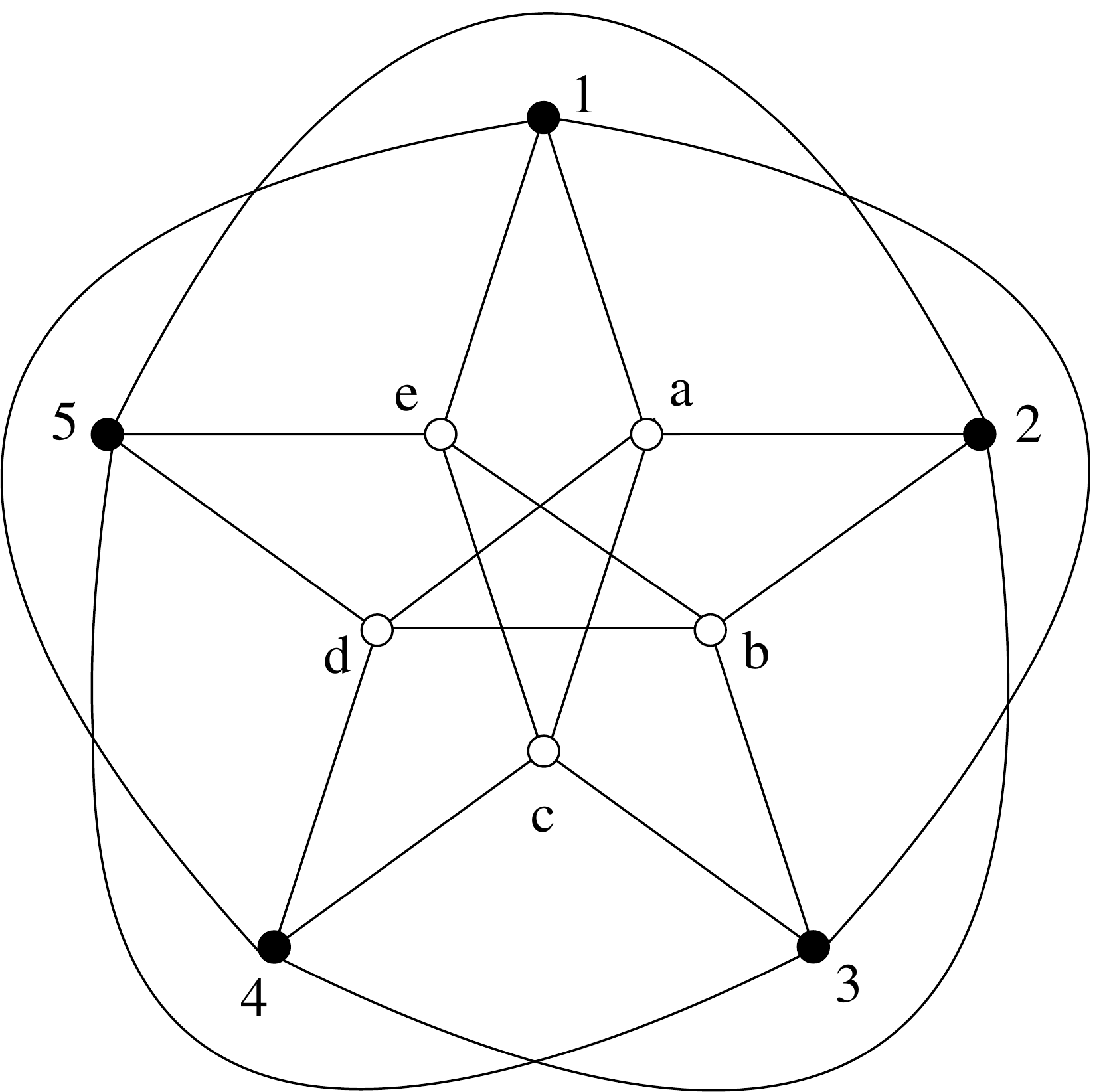}}
\end{center}
\caption{ $C_5^{--+}$}\label{C--+}
\end{figure}
\begin{figure}[ht]
\begin{center}
\scalebox{0.4}[.4]{\includegraphics{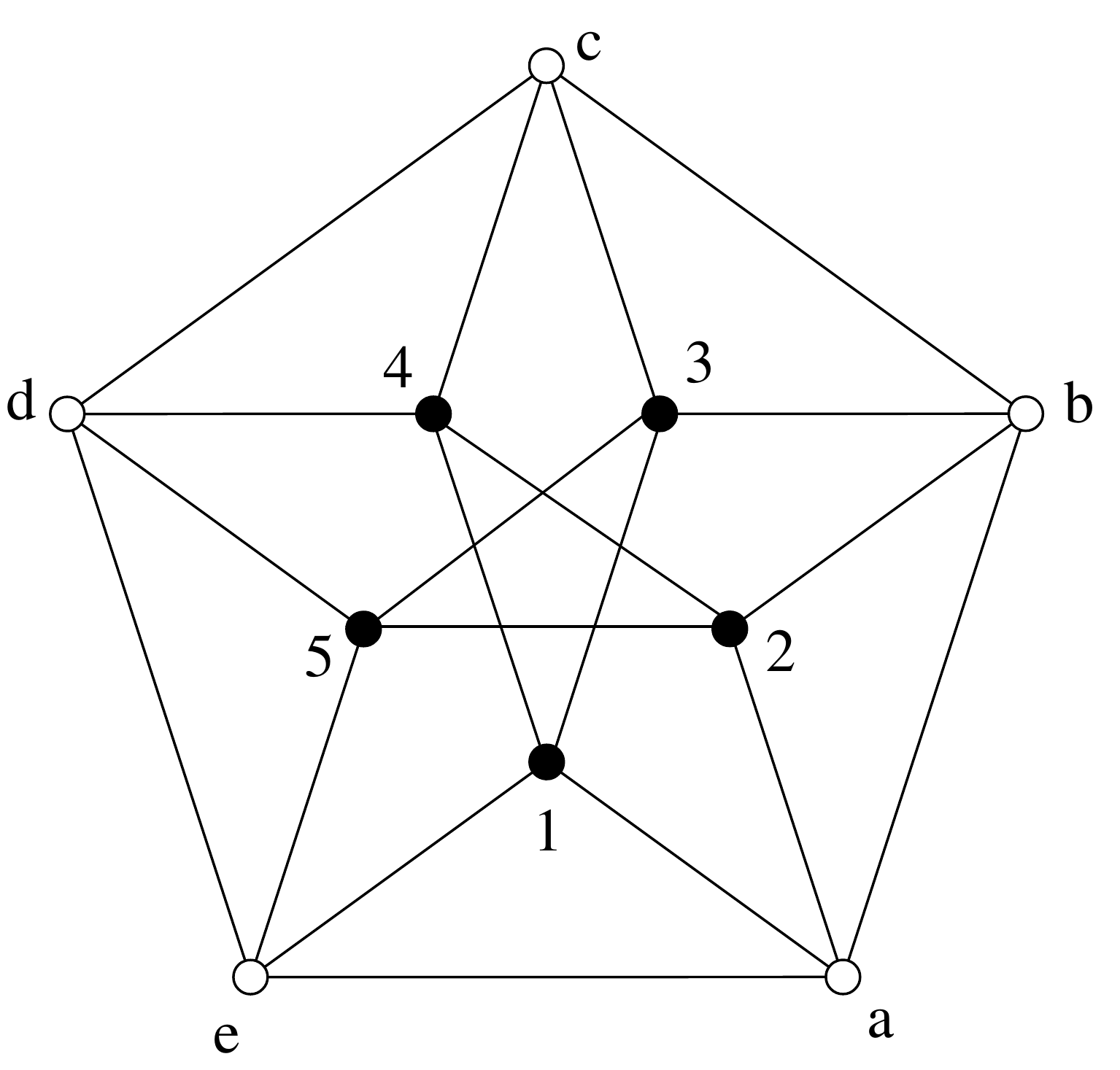}} \ \
\scalebox{0.4}[.4]{\includegraphics{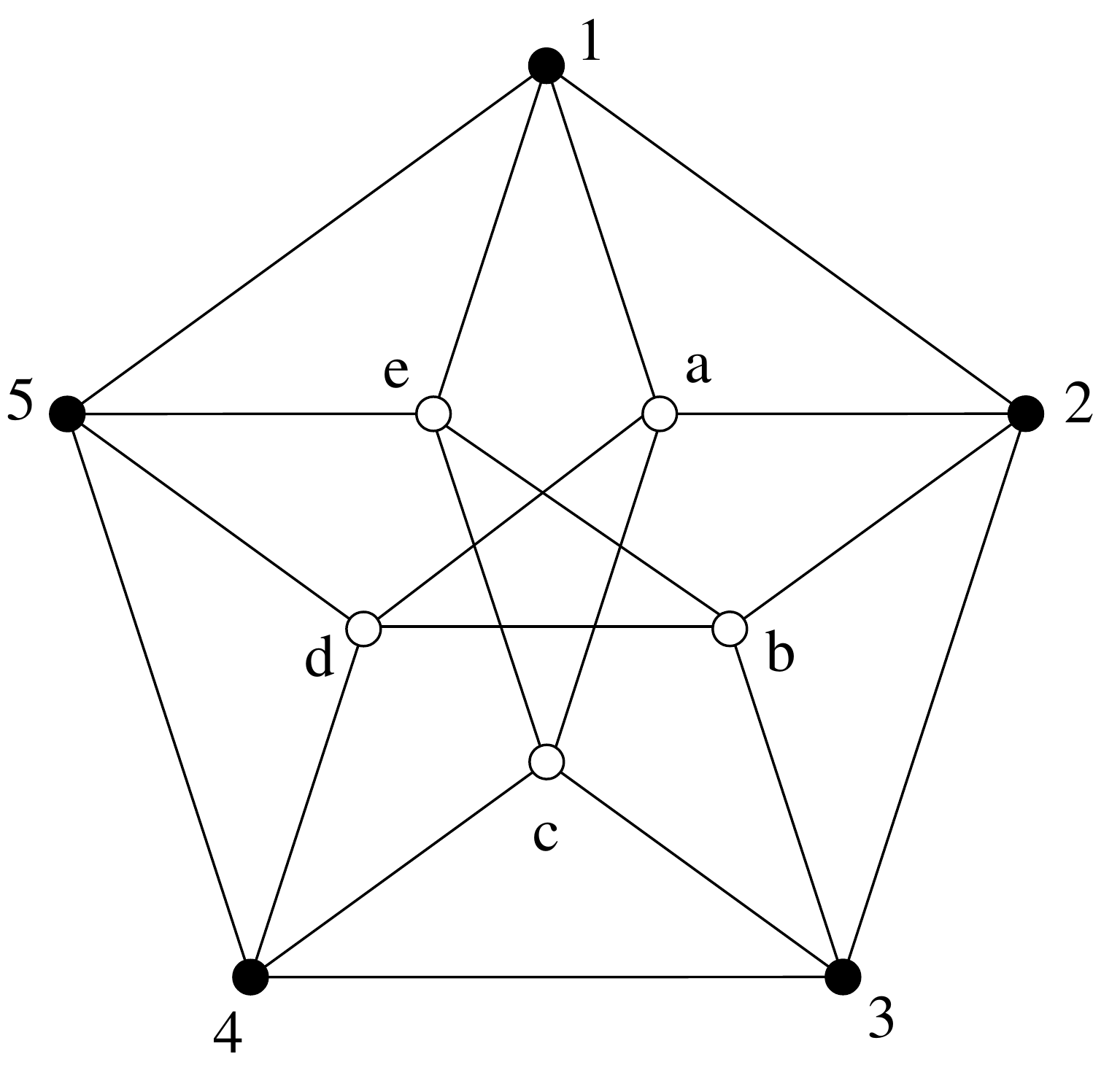}}
\end{center}
\caption{$C_5^{-++}$ (left) and $C_5^{+-+}$ (right)}\label{C-++/+-+}
\end{figure}
\\[1.5ex]
\indent
It is  easy to prove that $B(C_4)$ and $B^c(C_4)$ are isomorphic. Therefore $C_4^{xy+}$ and $C_4^{xy-}$ are isomorphic for any $x, y \in \{0, 1, +, -\}$.
\\[2ex]
\indent
It is also interesting to consider some transformations of the 5-cycle $C_5$ (the pentagon) because $C_5$ is isomorphic to its complement and  $C_5$ is also isomorphic to its line graph.
By the above formula for $C_n$, we have:
\\[1ex]
\indent
~~~~~~~~~~~~~$Sp(C_5) =
\{[\frac{1}{2}(5 + \sqrt{5})]^{(2)},~
[\frac{1}{2}(5 - \sqrt{5})]^{(2)},~ 0\}$.
\\[0.5ex]
Then
\\[1ex]
${\bf (p1)}$
$C_5^{+++}$ is a 4-regular graph (see Fig. \ref{C+++}) and
by Theorem \ref{th+++},
\\[1.5ex]
\indent
$Sp(C_5^{+++}) =
\{[\frac{1}{2}\langle 9 + \sqrt{5} + \sqrt{6 - 2\sqrt{5}}\rangle ]
^{(2)}, ~
[\frac{1}{2}\langle 9 + \sqrt{5} - \sqrt{6 - 2\sqrt{5}}\rangle  ]
^{(2)},~4,~
\\[1.5ex]
\indent
[\frac{1}{2}\langle 9 - \sqrt{5}  + \sqrt{6 - 2\sqrt{5}}\rangle  ]
^{(2)},~
[\frac{1}{2}\langle 9 - \sqrt{5}  - \sqrt{6 - 2\sqrt{5}}\rangle  ]
^{(2)},
~ 0\} $,
\\[2.5ex]
${\bf (p2)}$
$C_5^{--+}$ is a 4-regular graph (see Fig. \ref{C--+}) and
by Theorem \ref{th--+},
\\[1.5ex]
\indent
$Sp(C_5^{--+}) =
\{[\frac{1}{2}\langle 13 + \sqrt{5} + \sqrt {6 - 2\sqrt{5}}\rangle ]
^{(2)}, ~
[\frac{1}{2}\langle 13 + \sqrt{5} - \sqrt {6  - 2\sqrt{5}}\rangle ]
^{(2)},~
\\[1.5ex]
\indent
[\frac{1}{2}\langle 13- \sqrt{5}  + \sqrt {6 + 2\sqrt{5}}\rangle ]
^{(2)},~
[\frac{1}{2}\langle 13 - \sqrt{5}  - \sqrt {6 + 2\sqrt{5}}\rangle ]
^{(2)},~
4,~ 0\} $,
\\[2.5ex]
${\bf (p3)}$
$C_5^{-++}$ and $C_5^{+-+}$ are isomorphic  4-regular graphs (see Fig. \ref{C-++/+-+}) and by
Theorem \ref{th+-+} or \ref{th-++},
\\[1.6ex]
\indent
$Sp(C_5^{-++}) = Sp(C_5^{+-+}) =
 \{[\frac{1}{2}\langle 9  + \sqrt {11 + 2 \sqrt{5}}
 \rangle ]
^{(2)}, ~
[\frac{1}{2}\langle 9 + \sqrt {11 - 2 \sqrt{5}}\rangle  ]
^{(2)},~
\\[1.5ex]
\indent
4,~
[\frac{1}{2}\langle 9   + \sqrt{11 - 2\sqrt{5}}\rangle  ]
^{(2)},~
[\frac{1}{2}\langle 9   - \sqrt {11 + \sqrt{5}}\rangle  ]
^{(2)},
~ 0\} $.
\\[3ex]
\indent
By (c3) and (p3), $G^{-++}$ and $G^{+-+}$ are isomorphic if $G$ is either $C_4$ or $C_5$.
As we will see below,  a more general claim is true not only for
$C_4$ and $C_5$ but for any cycle $C$.
\begin{Theorem}
Let $G$ be an $r$-regular graph with $n$ vertices and $m$ edges. If $m=n$, then
$G^{xyz}$ and $G^{yxz}$ are isomorphic
for all $x, y, z \in \{ 0, 1, +, -\}$.
\end{Theorem}

\bp
By the Reciprocity Theorem \ref{reciprocity},
it is sufficient to prove our claim for $x, y \in \{0, 1, +, - \}$ and $z\in \{0, +\}$.
Since $m = n$ and  $n r = 2 m$, we have $r = 2$, and so  $G$ is 2-regular. Then $G$ is a disjoint union of  cycles.
If $A$ and $B$ are disjoint graphs,
then $(A\cup B)^{xy0} = A^{xy0} \cup A^{xy0}$ and $(A\cup B)^{xy+} = A^{xy+} \cup A^{xy+}$.
Therefore it is sufficient to prove our claim for a connected graph $G$. In this case $G$ is a cycle on
 $n$ vertices and  we can assume that  $V(G) = V = \{ v_1,  \cdots, v_n\}$ and
$E(G) = E = \{e_i: i= 1, \cdots, n\}$, where
$e_i = v_i v_{i+1}$ for  $i= 1, \cdots, n$ and $i+1$ is considered $\bmod~n$, and so $e_n = v_nv_1$.
Let for each $i$, $\alpha (v_i) = \sigma (v_{i+1}) = e_i$.
Then both $\alpha $ and $\sigma $ are isomorphisms from $G$ to $G^l$.
Recall  that $G^+ = G$, $G^- = G^c$, $G^0$ is the the  graph with $V(G^0) = V$ and with no edges , and   $G^1$ the a complete graph with $V(G^1) = V$.
Hence for every $x \in \{0, 1, +, -\}$, both $\alpha $ and $\sigma $ are isomorphisms from $G^x$ to $(G^l)^x$. Put
$\pi|_V = \alpha  \ \mbox{\ and \ } \ \pi|_E = \sigma ^{-1}$.
Since $G^{xy0}$ is a disjoint union of $G^x$ and $(G^l)^y$, we have: $\pi$ is an isomorphism
 from $G^{xy0}$ to $G^{yx0}$.

Now we show that $\pi$ is also an isomorphism from $G^{xy+}$ to $G^{yx+}$.
By definition of $G^{xy+}$,
$E(G^{xy+}) = E(G^x) \cup E((G^l)^y) \cup E(W)$, where
$E(W) = \{v_i e_i, v_{i+1}e_i: i = 1, \ldots , n\}$.
Recall that each $e_i = v_i v_{i+1}$.
Since $\pi(v_i) = \alpha(v_i) = e_i$ and
$\pi (e_j) = \sigma ^{-1} (e_j) = v_{j+1}$,
vertices $\pi(v_i)$ and $\pi(e_j)$ are adjacent in
$G^{yx+}$ if and only if $j+1 = i$ or $j+1 = i+1$
which is equivalent to $i= j+1$ or $i = j$.
Therefore $v_i$ and $e_j$ are adjacent in $G^{xy+}$ if and only if $\pi(v_i)$ and $\pi(e_j)$ are adjacent in $G^{yx+}$.
\ep

\section{Some remarks and questions}
\label{remarks}

\indent

${\bf (R1)}$
Each  factor of the Laplacian polynomials of $G^{xyz}$ ($x,y,z\in \{0,1,+,-\}$)  is
a polynomial in $\lambda $ of degree one or two. Therefore the explicit formula for the Laplacian spectrum and the number of spanning trees of
$G^{xyz}$ can be given in terms of those of $G$, respectively, as in
Corollaries \ref{Sp(00+)},
\ref{Sp(+0+)}, and \ref{Sp(0++)}.
\\[1.5ex]
\indent
${\bf (R2)}$
Let ${\cal R}$ be the set of regular graphs.
Obviously, if $G \in {\cal R}$, then $G^c \in {\cal R}$ and $G^l \in {\cal R}$.
If $G$ is an $r$-regular graph, then
$G^{+++}$ is $2r$-regular and $G^{---}$ is $(v(G) + e(G) - 2r -1)$-regular,
and so if $G \in {\cal R}$, then $G^{+++} \in {\cal R}$.
In other words, the set ${\cal R}$ of regular graphs is closed under
{\em $c$-operation}, {\em $l$-operation},
{\em $(+++)$-operation},
and {\em $(---)$-operation}.
Therefore using the corresponding results described above, one can give an algorithm (and the computer program) that for
any series $Z$ of $c$-, $l$-, $(+++)$-, and $(---)$-operations and the Laplacian spectrum $Sp(G)$ of any
$r$-regular graph $G$ provides the formula of the Laplacian spectrum of graph $F$ obtained from $G$ by the operation series $Z$ in terms of $r$, $v(G)$, and $Sp(G)$.
\\[1.5ex]
\indent
${\bf (R3)}$
Examples and results in Section \ref{examples} show that there exists a  regular graph $G$ such that  $G^{xyz}$ and
$G^{x'y'z'}$ are isomorphic although
$(x, y, z) \ne (x', y', z')$, where $x, y, z \in \{0, 1, +, -\}$.
It is also easy to see that
if $K$ is a complete graph, then
$K^{0yz} =  K^{-yz}$ and  $K^{x0z} = K^{x-z}$ as well as
$K^{1yz} =  K^{+yz}$ and  $K^{x1z} = K^{x+z}$.
\\[1.5ex]
\indent
${\bf (R4)}$
Suppose that a regular graph $G$ is uniquely defined by its Laplacian spectrum. Does it necessarily mean that $G^{xyz}$ is also uniquely defined by its Laplacian spectrum for every (or for some)
$x,y,z \in \{ +, -\}$ ?
\\[1.5ex]
\indent
\noindent{\bf Acknowledgement}:

Aiping Deng wishes to thank Michel Deza, Sergey V. Savchenko and Yaokun Wu for their kind help to bring about the cooperation with Alexander Kelmans.
We are thankful to the referees for useful remarks.

\newpage

\begin{center}
{\LARGE Appendix}
\end{center}

Let $G$ be an $r$-regular graph with $n$ vertices and $m$ edges and let $s = n + m$,
and so $2m = rn$ and $s$ is the number of vertices of $G^{xyz}$.
The tables below provide the formulas for
$L(\lambda, G^{xyz})$ and $t(G^{xyz})$ for all $x, y, z \in \{0, 1, + , -\}$ in terms of $n$, $m$, $r$, and the Laplacian eigenvalues of $G$.
Obviously, $t(G^{xyz}) = 0$ if and only if
$G^{xyz}$ is not connected.
Note that for $z = 0$,  graph $G^{xyz}$ is not connected, and so $t(G^{xyz}) = 0$. For $z=1$, graph $G^{xyz}$ is  connected, and so
$t(G^{xyz}) > 0$.
\\[3ex]
\indent
{\sc The list of }
$L(\lambda, G^{xyz})$ {\sc and} $t(G^{xyz})$
{\sc with}
$z=0$.
\\[1.5ex]
\begin{tabular}{|l|l|l|l|}
  \hline
  & $xyz$ &  $L(\lambda, G^{xyz})$ & $t(G^{xyz})$  \\
 \hline
1 & $0~0~0$ & $\lambda^{m+n}$ & ~~~0 \\
 \hline
2 & $1~0~0$ & $\lambda^{m+1}(\lambda-n)^{n-1}$ & ~~~0  \\
 \hline
3 & $+0~0$ & $\lambda^m L(\lambda,G)$ & ~~~0  \\
 \hline
4 & $-0~0$ & $(-1)^n
(\lambda - n)^{-1} \lambda^{m+1} L(n-\lambda, G)$ & ~~~0 \\
 \hline
5 & $0~1~0$ & $\lambda^{n+1} (\lambda - m)^{m-1} $ & ~~~0  \\
 \hline
6 & $1~1~0$ & $\lambda^2 (\lambda - m)^{m-1} (\lambda -n)^{n-1} $ & ~~~0  \\
 \hline
7 & $+1~0$ & $\lambda (\lambda - m)^{m-1} L(\lambda, G)$ & ~~~0   \\
 \hline
8 & $-1~0$  & $(-1)^n \lambda^2 (\lambda -n)^{-1} (\lambda -m)^{m-1} L(n-\lambda, G)$ & ~~~0  \\
 \hline
9 & $0+0$  & $\lambda ^n (\lambda - 2r)^{m - n} L(\lambda , G)$ &~~~0  \\
 \hline
10 & $1+0$  & $\lambda (\lambda - n)^{n-1}(\lambda - 2r)^{m - n} L(\lambda , G)$ & ~~~0  \\
 \hline
11 & $++0$  & $(\lambda - 2r)^{m - n} L(\lambda , G)^2$ & ~~~0  \\
 \hline
12 & $-+0$  & $(-1)^n (\lambda - 2r)^{m - n}
\lambda (\lambda - n)^{-1}
L(n - \lambda , G) L(\lambda , G)$ &  ~~~0 \\
 \hline
13 & $0-0$  & $(-1)^n (\lambda - m)^{-1} \lambda^{n+1} (\lambda - m + 2r)^{m-n}
L(m-\lambda, G)$ &  ~~~0 \\
 \hline
14 & $1-0$  & $ (-1)^n (\lambda - m)^{-1} \lambda^2 (\lambda -n)^{n-1} (\lambda - m + 2r)^{m-n} L(m-\lambda, G)$ &  ~~~0 \\
 \hline
15 & $+-0$  & $(-1)^n
\lambda (\lambda - m)^{-1}
(\lambda - m + 2r)^{m - n} L(m - \lambda , G)
L(\lambda , G)$ &  ~~~0  \\
 \hline
 16 & $--0$ & $\lambda ^2( \lambda - m)^{-1}(\lambda - n)^{-1}
(\lambda - m +  2r)^{m - n}
L(m - \lambda , G) L(n - \lambda , G)$ & ~~~0  \\
 \hline
\end{tabular}
\\[3ex]

\newpage
{\sc The list of }
$L(\lambda, G^{xyz})$ {\sc and} $t(G^{xyz})$
{\sc with}
$z=1$.
\\[1.5ex]
\begin{tabular}{|l|l|l|l|}
  \hline
  & $xyz$ &  $L(\lambda, G^{xyz})$ & $t(G^{xyz})$  \\
 \hline
1 & $0~0~1$ & $\lambda (\lambda - s)(\lambda - n)^{m-1} (\lambda - m)^{n-1}$ & $n^{m-1} m^{n-1}$ \\
 \hline
2 & $1~0~1$ & $\lambda (\lambda - s)^n (\lambda - n)^{m-1} $ & $n^{m-1} s^{n-1}$   \\
 \hline
3 & $+0~1$ & $(\lambda - m)^{-1}
(\lambda - s)(\lambda -n)^{m-1} L(\lambda - m, G)$ & $n^{m-1} \prod_{i=1}^{n-1}(m+ \lambda_i)$  \\
 &  &  $=\lambda (\lambda -s)(\lambda - n)^{m-1}\prod_{i=1}^{n-1}(\lambda - m - \lambda_i)$ &  \\
 \hline
4 & $-0~1$ & $(-1)^{n}\lambda (\lambda - n)^{m-1} L(s -\lambda, G) $ &  $n^{m-1}
\prod_{i=1}^{n-1}(s - \lambda_i)$ \\
 &  &  $=\lambda (\lambda - s)(\lambda - n)^{m-1} \prod_{i=1}^{n-1}(\lambda - s + \lambda_i)$ & \\
 \hline
5 & $0~1~1$ & $\lambda (\lambda - s)^{m} (\lambda - m)^{n -1} $ & $s^{m-1}m^{n-1}$  \\
 \hline
6 & $1~1~1$ & $\lambda (\lambda - s)^{s-1} $ & $s^{s -2}$  \\
 \hline
7 & $+1~1$ & $\lambda (\lambda - m)^{-1}(\lambda - s)^m L(\lambda - m, G)$ & $s^{m -1}\prod_{i=1}^{n-1}(m + \lambda_i)$ \\
 &  & $= \lambda (\lambda - s)^m \prod_{i=1}^{n-1}(\lambda - m - \lambda_i)$ &  \\
 \hline
8 & $-1~1$  & $(-1)^n \lambda (\lambda - s)^{m-1} L(s - \lambda, G)$ & $s^{m-1} \prod_{i=1}^{n-1}(s - \lambda_i)$  \\
&  &  $= \lambda (\lambda - s)^m \prod_{i=1}^{n-1}(\lambda - s + \lambda_i)$  &   \\
 \hline
9 & $0+1$  & $\lambda (\lambda - s) (\lambda - n)^{-1} (\lambda - m)^{n-1}$ & $m^{n-1}
(n+ 2r)^{m-n}\prod_{i=1}^{n-1} (n + \lambda_i)$  \\
 &  &  $(\lambda - n - 2r)^{m-n}
L(\lambda - n, G)$  & $ $ \\
 \hline
10 & $1+1$  & $\lambda (\lambda - n)^{-1} (\lambda - s)^n (\lambda - n - 2r)^{m-n}$ &
$s^{n-1} (n + 2r)^{m-n}$  \\
& & $L(\lambda - n, G)$ & $\prod_{i=1}^{n-1}(n + \lambda_i)$   \\
 \hline
11 & $++1$  & $\lambda (\lambda - s) (\lambda -m)^{-1} (\lambda - n)^{-1}$ & $(n+ 2r)^{m-n}$  \\
 &  &  $(\lambda - n - 2r)^{m-n} L(\lambda - m, G) L(\lambda - n, G)$ & $\prod_{i+1}^{n-1}(m + \lambda_i)(n + \lambda_i)$ \\
 \hline
12 & $-+1$  & $\lambda (\lambda - n)^{-1} (\lambda - n - 2r)^{m-n}$ &  $(n + 2r)^{m-n}$ \\
 & & $ L(\lambda - n, G) L(s - \lambda, G)$ &  $\prod_{i = 1}^{n-1}(n+ \lambda_i)(s - \lambda_i)$ \\
 \hline
13 & $0-1$  & $(-1)^{n}\lambda (\lambda - m)^{n-1} (\lambda - s + 2r)^{m-n}L(s - \lambda, G) $
&  $m^{n-1}(s -2r)^{m-n}\prod_{i=1}^{n-1}(s - \lambda_i)$  \\
 \hline
14 & $1-1$  & $(-1)^n \lambda (\lambda -s)^{n-1}(\lambda - s + 2r)^{m-1}L(s - \lambda, G)$
& $s^{n-1}(s -2r)^{m-n}\prod_{i=1}^{n-1}(s - \lambda_i)$ \\
 \hline
 15 & $+-1$  & $(-1)^n \lambda (\lambda-m)^{-1} (\lambda - s + 2r)^{m-n}$ & $(s -2r)^{m-n} $  \\
& & $L(s - \lambda, G) L(\lambda - m, G)$ & $\prod_{i=1}^{n-1}(s - \lambda_i)(m + \lambda_i) $  \\
 \hline
16 & $--1$ & $\lambda (\lambda - s)^{-1}
(\lambda - s + 2r)^{m-n}
L(s - \lambda, G)^2$ & $(s - 2r)^{m-n}
\prod_{i=1}^{n-1}(s - \lambda_i)^2$  \\
 \hline
\end{tabular}

\newpage
{\sc The list of }
$L(\lambda, G^{xyz})$ {\sc and } $t(G^{xyz})$
{\sc with }
$z=+$.
\\[1.5ex]
\begin{tabular}{|l|l|l|l|}
  \hline
  & $xyz$ &  $L(\lambda, G^{xyz})$ & $t(G^{xyz})$  \\
 \hline
1 & $0~0+$ & $\lambda (\lambda - r -2)(\lambda -2)^{m-n} \prod_{i=1}^{n-1}\{\lambda^2 - \lambda(r+2) + \lambda_i\}$
& $ns^{-1}(r + 2) 2^{m - n} t(G)$ \\
 \hline
2 & $1~0+$ & $ \lambda (\lambda  - r - 2)(\lambda  - 2)^{m-n}$ &
$s^{-1}(r + 2) 2^{m - n}
\prod_{i=1}^{n-1}(2n + \lambda_i)$  \\
 & &
 $
 \prod_{i=1}^{n-1}
 \{(\lambda -2)(\lambda - n -r) -2r + \lambda_i \}$ & \\
 \hline
3 & $+0+$ & $\lambda (\lambda  - r - 2) (\lambda  - 2)^{m-n}$ &
$ns^{-1}(r + 2) 2^{m - n} 3^{n-1} t(G)$ \\
&  &  $\prod_{i=1}^{n-1} \{ (\lambda - 2) (\lambda - r - \lambda_i) - 2r + \lambda_i \}$  &  \\
 \hline
4 & $-0+$ & $\lambda (\lambda - r -2)(\lambda - 2)^{m-n}$  &
$s^{-1}(r + 2) 2^{m - n}
\prod_{i=1}^{n-1}(2n - \lambda_i)$  \\
&  &  $\prod_{i=1}^{n-1} \{(\lambda - 2)(\lambda - n -r + \lambda_i) - 2r + \lambda_i\}$  &  \\
 \hline
5 & $0~1+$ & $ \lambda (\lambda  - r - 2)(\lambda - m - 2)^{m-n}$ &
$s^{-1}(r + 2) (m + 2)^{m - n}$  \\
 &  &  $ \prod_{i=1}^{n-1} (\lambda - r)(\lambda -m -2) -2r
 + \lambda _i \} $  &  $\prod_{i=1}^{n-1}(mr + \lambda_i)$ \\
 \hline
6 & $1~1+$ & $ \lambda (\lambda  - r - 2)(\lambda - m - 2)^{m-n}$ &
$s^{-1}(r + 2) (m + 2)^{m - n}$  \\
 &  &  $ \prod_{i=1}^{n-1} \{(\lambda - n -r)(\lambda - m -2) -2r + \lambda_i \}$  & $\prod_{i=1}^{n-1}(mn + 2n + mr + \lambda_i)$ \\
 \hline
7 & $+1+$ & $\lambda (\lambda  - r - 2)(\lambda - m - 2)^{m-n}$ &
$s^{-1}(r + 2) (m + 2)^{m - n}$   \\
 &  &  $ \prod_{i=1}^{n-1} \{(\lambda - r - \lambda_i)(\lambda - m -2) -2r + \lambda_i \}$  & $\prod_{i=1}^{n-1}(mr + m\lambda_i + 3\lambda_i)$ \\
 \hline
8 & $-1+$  & $\lambda (\lambda  - r - 2)(\lambda - m - 2)^{m-n}$ &
$s^{-1}(r + 2) (m + 2)^{m - n}$  \\
 &   &   $ \prod_{i=1}^{n-1} \{(\lambda - n -r + \lambda_i)(\lambda - m -2) -2r + \lambda_i \}$
 & $\prod_{i=1}^{n-1}(mr + mn  + 2n - m\lambda_i - \lambda_i)$  \\
 \hline
9 & $0++$  & $\lambda
 (\lambda - r- 2) (\lambda - 2r - 2)^{m-n}$ &
 $ns^{-1} (r+2) 2^{m-n} (r+1)^{m-1}t(G)$  \\
& & $\prod_{i=1}^{n-1}\{ (\lambda - r)(\lambda - 2 - \lambda _i) - 2r + \lambda _i \}$ &   \\
 \hline
10 & $1++$  & $ \lambda (\lambda  - r - 2)(\lambda  - 2r - 2)^{m-n}$ &
$s^{-1}(r + 2) (2r + 2)^{m - n}$  \\
& &
 $\prod_{i=1}^{n-1} \{(\lambda - n -r)(\lambda -2 - \lambda_i) -2r + \lambda_i \}$ & $\prod_{i=1}^{n-1}\{2n + \lambda_i(n + r + 1)\}$   \\
 \hline
11 & $+++$  & $\lambda (\lambda  - r - 2)(\lambda - 2r -2)^{m - n}
$ &
$ns^{-1}(r+2)(2r+2)^{m-n} t(G)$  \\
& & $\prod_{i=1}^{n - 1}
\{(\lambda - r - \lambda _i)(\lambda - 2 - \lambda _i) - 2r + \lambda _i\}$ & $\prod_{i=1}^{n-1}(r + 3 + \lambda_i)$  \\
 \hline
12 & $-++$  & $\lambda(\lambda-r-2)(\lambda -2r-2)^{m-n}$ &
$s^{-1}(r+2)(2r+2)^{m-n}$  \\
& & $\prod_{i=1}^{n-1} \{(\lambda -n -r +\lambda_i)(\lambda  - 2  - \lambda_i)- 2r + \lambda_i\}$ & $\prod_{i=1}^{i-1}(2n+ \lambda_i(n + r -1 - \lambda_i))$  \\
 \hline
13 & $0-+$  & $\lambda (\lambda - r -2) (\lambda - m + 2r - 2)^{m-n}$ &
$s^{-1}(r+2)(m-2r+2)^{m-n}$  \\
& & $\prod_{i=1}^{n-1} \{(\lambda - r) (\lambda - m  -2 + \lambda_i)- 2r + \lambda_i\}$ & $\prod_{i=1}^{n-1}( mr + \lambda_i r + \lambda_i)$  \\
 \hline
14 & $1-+$  & $ \lambda (\lambda  - r - 2)(\lambda -m + 2r - 2)^{m-n}$ &
$s^{-1}(r+2)(m-2r+2)^{m-n}$ \\
& & $\prod_{i=1}^{n-1} \{(\lambda - n -r)(\lambda - m -2 + \lambda_i) -2r + \lambda_i \}$ & $\prod_{i=1}^{n-1}\{ mn + 2n + mr - \lambda_i(n +r -1)\}$  \\
 \hline
15 & $+-+$  & $\lambda (\lambda - r - 2)
(\lambda -m +2r -2)^{m - n}
$ &
$s^{-1} (r+2)(m-2r+2)^{m-n}$  \\
& & $\prod_{i=1}^{n-1}\{ (\lambda - r - \lambda_i) (\lambda - m - 2 + \lambda _i) - 2r + \lambda_i\}$
&  $\prod_{i=1}^{n-1}\{mr + \lambda_i(m - r + 3 - \lambda_i)\}$  \\
 \hline
16 & $--+$  & $\lambda (\lambda  - r - 2) (\lambda  - m + 2r - 2)^{m - n}
\prod_{i=1}^{n-1}$ &
$s^{-1}(r+2)(m-2r+2)^{m-n}
\prod_{i=1}^{n-1}\{mn $ \\
& & $ \{(\lambda - n - r + \lambda_i) (\lambda - m - 2 + \lambda _i) - 2r + \lambda_i\}$
& $+ 2n + mr + \lambda_i(\lambda_i - m - n -r -1)\}$  \\
 \hline
\end{tabular}

\newpage

{\sc The list of }
$L(\lambda, G^{xyz})$ {\sc and } $t(G^{xyz})$
{\sc with }
$z=-$.
\\[1.5ex]
\begin{tabular}{|l|l|l|l|}
  \hline
  & $xyz$ &  $L(\lambda, G^{xyz})$ & $t(G^{xyz})$  \\
 \hline
1 & $0~0-$ & $\lambda (\lambda - s + r + 2)(\lambda - n + 2)^{m-n}$ & $s^{-1}(s -r -2)(n -2)^{m-n}$  \\
 &   &   $\prod_{i=1}^{n-1}\{ (\lambda - m + r)(\lambda -n + 2) -2r + \lambda_i\}$  & $\prod_{i=1}^{n-1}(mn - nr - 2m + \lambda_i)$ \\
 \hline
2 & $1~0-$ & $\lambda (\lambda - s + r + 2)(\lambda - n + 2)^{m-n}$ & $s^{-1}(s -r -2)(n -2)^{m-n}$  \\
 &   &   $\prod_{i=1}^{n-1}\{ (\lambda - s + r)(\lambda -n + 2) -2r + \lambda_i\}$  & $\prod_{i=1}^{n-1}(ns - nr - 2s + \lambda_i)$ \\
 \hline
3 & $+0-$ & $\lambda (\lambda - s + r + 2)(\lambda - n + 2)^{m-n}$ & $s^{-1}(s -r -2)(n -2)^{m-n}$  \\
 &   &   $\prod_{i=1}^{n-1}\{ (\lambda - m + r - \lambda_i)(\lambda -n + 2) -2r + \lambda_i\}$  & $\prod_{i=1}^{n-1}(mn - nr - 2m + n\lambda_i - \lambda_i)$ \\
 \hline
4 & $-0-$ & $\lambda (\lambda - s + r + 2)(\lambda - n + 2)^{m-n}$ & $s^{-1}(s -r -2)(n -2)^{m-n}$  \\
 &   &   $\prod_{i=1}^{n-1}\{ (\lambda - s + r + \lambda_i)(\lambda -n + 2) -2r + \lambda_i\}$  & $\prod_{i=1}^{n-1}(ns - nr - 2s -n\lambda_i + 3\lambda_i)$ \\
 \hline
5 & $0~1-$ & $\lambda (\lambda - s + r +2)(\lambda - s + 2)^{m-n}$ & $s^{-1}(s -r -2)(s -2)^{m-n}$  \\
 &  &  $\prod_{i=1}^{n-1}\{(\lambda -s + 2)(\lambda -m + r) -2r + \lambda_i\}$  & $\prod_{i=1}^{n-1}(ms - rs - 2m + \lambda_i)$ \\
 \hline
6 & $1~1-$ & $\lambda (\lambda - s + r +2)(\lambda - s + 2)^{m-n}$ & $s^{-1}(s -r -2)(s -2)^{m-n}$ \\
 &  &  $\prod_{i=1}^{n-1}\{(\lambda -s)(\lambda -s + r + 2) + \lambda_i\}$  & $\prod_{i=1}^{n-1}(s(s- r - 2) + \lambda_i)$ \\
 \hline
7 & $+1-$ & $\lambda (\lambda - s + r +2)(\lambda - s + 2)^{m-n}$ & $s^{-1}(s -r -2)(s -2)^{m-n}$  \\
 &  &  $\prod_{i=1}^{n-1}\{(\lambda -s + 2)(\lambda -m + r - \lambda_i) -2r + \lambda_i\}$  & $\prod_{i=1}^{n-1}(s(m - r + \lambda_i) - 2m - \lambda_i)$ \\
 \hline
8 & $-1-$  & $\lambda (\lambda - s + r +2)(\lambda - s + 2)^{m-n}$ & $s^{-1}(s -r -2)(s -2)^{m-n}$  \\
 &  &  $\prod_{i=1}^{n-1}\{(\lambda -s + 2)(\lambda -s + r + \lambda_i) -2r + 3\lambda_i\}$
 & $\prod_{i=1}^{n-1}(s(s - r - \lambda_i -2) + \lambda_i)$ \\
 \hline
9 & $0+-$  & $\lambda (\lambda - s + r + 2)(\lambda - n - 2r + 2)^{m-n}$ & $s^{-1}(s -r -2)(n + 2r -2)^{m-n}$  \\
 &  &  $\prod_{i=1}^{n-1}\{ (\lambda - m + r)(\lambda - n + 2 - \lambda _i) - 2r + \lambda_i \}$
 & $\prod_{i=1}^{n-1}((m - r)(n + \lambda_i) - 2m + \lambda_i)$ \\
 \hline
10 & $1+-$  & $\lambda (\lambda - s + r + 2)(\lambda - n - 2r + 2)^{m-n}$ & $s^{-1}(s -r -2)(n + 2r -2)^{m-n}$  \\
 &  &  $\prod_{i=1}^{n-1}\{ (\lambda - s + r)(\lambda - n + 2 - \lambda _i) - 2r + \lambda_i \}$
 & $\prod_{i=1}^{n-1}((s - r)(n + \lambda_i) - 2s + \lambda_i)$ \\
 \hline
11 & $++-$  & $\lambda (\lambda - s + r + 2)(\lambda - n - 2r + 2)^{m-n}$ & $s^{-1}(s -r -2)(n + 2r -2)^{m-n}$  \\
 &  &  $\prod_{i=1}^{n-1}\{ (\lambda - m + r - \lambda_i)(\lambda - n + 2 - \lambda _i) - 2r + \lambda_i \}$
 & $\prod_{i=1}^{n-1}((m - r + \lambda_i)(n + \lambda_i) - 2m - \lambda_i)$ \\
 \hline
12 & $-+-$ & $\lambda (\lambda - s + r + 2)(\lambda - n - 2r + 2)^{m-n}$ & $s^{-1}(s -r -2)(n + 2r -2)^{m-n}$  \\
 &  &  $\prod_{i=1}^{n-1}\{ (\lambda - s + r + \lambda_i)(\lambda - n + 2 - \lambda _i) - 2r + \lambda_i \}$
 & $\prod_{i=1}^{n-1}((s - r - \lambda_i)(n + \lambda_i) - 2s + 3\lambda_i)$ \\
 \hline
13 & $0--$  & $\lambda (\lambda - s + r + 2)(\lambda - s + 2r + 2)^{m-n}$ & $s^{-1}(s -r -2)(s - 2r -2)^{m-n}$ \\
& & $\prod_{i=1}^{n-1}\{ (\lambda - m + r) (\lambda - s + 2 + \lambda_i) - 2r + \lambda_i \}$
& $\prod_{i=1}^{n-1}((m - r)(s - \lambda_i) - 2m + \lambda_i)$  \\
 \hline
14 & $1--$  & $\lambda (\lambda - s + r + 2)(\lambda - s + 2r + 2)^{m-n}$ & $s^{-1}(s -r -2)(s - 2r -2)^{m-n}$ \\
& & $\prod_{i=1}^{n-1}\{ (\lambda - s + r) (\lambda - s + 2 + \lambda_i) - 2r + \lambda_i \}$
& $\prod_{i=1}^{n-1}((s - r)(s - \lambda_i) - 2s + \lambda_i)$   \\
 \hline
15 & $+--$  & $\lambda (\lambda - s + r + 2)(\lambda - s + 2r + 2)^{m-n}$ & $s^{-1}(s -r -2)(s - 2r -2)^{m-n}$ \\
& & $\prod_{i=1}^{n-1}\{ (\lambda - m + r - \lambda_i) (\lambda - s + 2 + \lambda_i) - 2r + \lambda_i \}$
& $\prod_{i=1}^{n-1}((m - r + \lambda_i)(s - \lambda_i) - 2m - \lambda_i)$   \\
 \hline
16 & $---$ & $\lambda (\lambda - s + r + 2)(\lambda - s + 2r + 2)^{m-n}$ & $s^{-1}(s -r -2)(s - 2r -2)^{m-n}$ \\
& & $\prod_{i=1}^{n-1}\{ (\lambda - s + r + \lambda_i) (\lambda - s + 2 + \lambda_i) - 2r + \lambda_i \}$
& $\prod_{i=1}^{n-1}((s - r - \lambda_i)(s - \lambda_i) - 2s + 3\lambda_i)$   \\
 \hline
\end{tabular}

\end{document}